\documentclass[a4paper]{scrartcl}

\usepackage[utf8]{inputenc}
\usepackage[T1]{fontenc}
\usepackage{lmodern}
\usepackage[svgnames,dvipsnames,rgb]{xcolor} 
\usepackage{amsfonts}
\usepackage{amsmath} 
\usepackage{amssymb}
\usepackage{amsthm}
\usepackage{graphicx}
\usepackage{hyperref}
\hypersetup{linktocpage=true,colorlinks=true,linkcolor=RoyalBlue,citecolor=PineGreen,urlcolor=RoyalBlue}
\usepackage{float}
\usepackage{tikz,pgf}
\usetikzlibrary{calc}
\usetikzlibrary{patterns}
\usetikzlibrary{angles}
\usetikzlibrary{backgrounds}

\usepackage{cleveref}
\usepackage{subcaption}
\usepackage{stackrel}

\usepackage[shortlabels]{enumitem}

\newtheorem{theorem}{Theorem}[section]
\newtheorem{lemma}[theorem]{Lemma}
\newtheorem{proposition}[theorem]{Proposition}
\newtheorem{corollary}[theorem]{Corollary}

\theoremstyle{definition}\newtheorem{definition}[theorem]{Definition}
\theoremstyle{definition}\newtheorem{example}[theorem]{Example}
\theoremstyle{definition}\newtheorem{remark}[theorem]{Remark}

\newcommand{\FlexRiLoG}{\textsc{FlexRiLoG}}

\newcommand{\oriented}[2]{(#1,#2)}

\newcommand{\blue}{\text{blue}}
\newcommand{\red}{\text{red}}

\DeclareMathOperator{\NAC}{NAC}
\newcommand{\nac}[1]{\NAC_{#1}}

\newcommand{\NN}{\mathbb{N}}
\newcommand{\RR}{\mathbb{R}}
\newcommand{\CC}{\mathbb{C}}
\newcommand{\QQ}{\mathbb{Q}}

\newcommand{\KK}{\mathbb{K}}

\newcommand{\C}{\mathcal{C}}

\newcommand{\Rot}{\Theta}

\newcommand{\Ck}{\mathcal{C}_k} 
\newcommand{\Aut}{\operatorname{Aut}}
\colorlet{ecol}{black!50!white}
\definecolor{colR}{rgb}{.932,.172,.172} 
\definecolor{colB}{rgb}{.255,.41,.884} 
\colorlet{colG}{Gold}
\colorlet{col1}{LightGreen}
\colorlet{col2}{IndianRed!80!white}
\colorlet{col3}{Gold}
\colorlet{col4}{LightSkyBlue}
\colorlet{col5}{BurlyWood}
\colorlet{col6}{Purple!60!white}
\colorlet{ncol}{DarkSeaGreen}

\tikzstyle{gvertex}=[circle, draw, fill=black, inner sep=0pt, minimum size=4pt]
\tikzstyle{smallgvertex}=[circle, draw, fill=black, inner sep=0pt, minimum size=2pt]
\tikzstyle{midgvertex}=[circle, draw, fill=black, inner sep=0pt, minimum size=3pt]
\tikzstyle{joint}=[circle, draw=ecol, fill=ecol, inner sep=0pt, minimum size=1pt]

\colorlet{colvR}{black!70!white}
\tikzstyle{fvertex}=[circle,thick,draw=colvR,fill=white,inner sep=0pt, minimum size=4.5pt]
\tikzstyle{midfvertex}=[circle,thick,draw=colvR, fill=white, inner sep=0pt, minimum size=3.25pt]
\tikzstyle{smallfvertex}=[circle,thick,draw=colvR, fill=white, inner sep=0pt, minimum size=2pt]

\tikzstyle{edge}=[line width=1.5pt,ecol]
\tikzstyle{redge}=[edge,colR]
\tikzstyle{bedge}=[edge,colB]
\tikzstyle{gedge}=[edge,colG]

\tikzstyle{gridl}=[ecol]
\tikzstyle{gridp}=[inner sep=1pt,circle,fill=black!70!white]
\tikzstyle{sym}=[ecol,dashed]

\tikzstyle{axes}=[gridl,-latex]

\tikzstyle{ngvertex}=[gvertex, draw=ncol, fill=ncol]
\tikzstyle{edgeq}=[edge,gray!60,densely dashed]
\tikzstyle{nedge}=[edge,ncol]
\tikzstyle{bdedge}=[line width=1.5pt,colB, densely dashed]
\tikzstyle{rdedge}=[line width=1.5pt,colR, densely dashed]

\tikzstyle{genericgraph}=[dashed]
\tikzstyle{labelsty}=[font=\scriptsize]
\tikzstyle{indicatededge}=[pin={[pin distance=6pt,pin edge={thin}]20:},pin={[pin distance=6pt,pin edge={thin}]-20:},pin={[pin distance=10pt,pin edge={thin}]2:}]

\tikzstyle{ribbon}=[densely dashed,rounded corners,line width=1pt,shorten <= -6pt, shorten >= -6pt]
\tikzstyle{brace}=[line width=1pt,ecol]

\tikzstyle{category}=[font=\small]
\tikzstyle{interiorangle}=[angle radius=0.3cm]

\tikzstyle{thin}=[col3]
\tikzstyle{fat}=[col1]

\title{Flexing infinite frameworks with applications to braced Penrose tilings}
\author{%
Sean Dewar$^{\text{{\tiny $\square$}}, \dagger}$
\and
Jan Legersk\'y$^{\circ, \diamond, \ast}$
}

\date{}

\begin{document}
\maketitle

\footnotetext{\hspace{0.15cm}$^{\text{{\tiny $\square$}}}$ Supported by the Austrian Science Fund (FWF): P31888.\\%
$^\dagger$ Supported by the Fields Institute for Research in Mathematical Sciences.\\%
$^\circ$ Supported by the Austrian Science Fund (FWF): P31061.\\%
$^\diamond$ Supported by the Ministry of Education, Youth and Sports of the Czech Republic, project no. CZ.02.1.01/0.0/0.0/16\_019/0000778.\\%
$^{\ast}$ Corresponding author (\texttt{jan.legersky@fit.cvut.cz})\\
This is the accepted version of the
following article: \textit{Sean Dewar and Jan Legersk\'y.
Flexing infinite frameworks with applications to braced Penrose tilings. Discrete Applied Mathematics. 324:1--17, 2023},
which has been published in final form at \href{https://doi.org/10.1016/j.dam.2022.09.002}{doi:10.1016/j.dam.2022.09.002}.\\
Licensed under Creative Commons License CC-BY-NC-ND.
}

\begin{abstract}
	A planar framework --- a graph together with a map of its vertices to the plane ---
	is flexible if it allows a continuous deformation preserving the distances between adjacent vertices.
	Extending a recent previous result,
	we prove that a connected graph with a countable vertex set can be realized as a flexible framework if and only if it has a so-called NAC-coloring.
	The tools developed to prove this result are then applied to frameworks where every 4-cycle is a parallelogram,
	and countably infinite graphs with $n$-fold rotational symmetry.
	With this, we determine a simple combinatorial characterization that determines whether the 1-skeleton of a Penrose rhombus tiling with a given set of braced rhombi will have a flexible motion, and also whether the motion will preserve 5-fold rotational symmetry.
\end{abstract}

\section{Introduction}

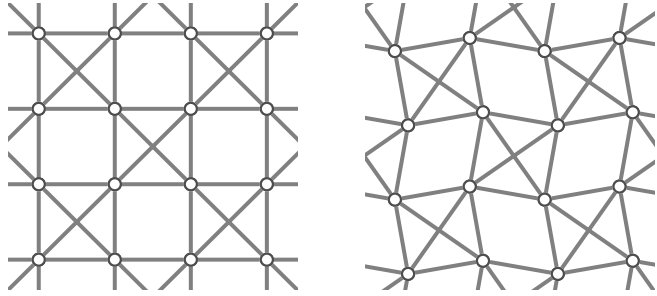
\begin{figure}[ht]
	\centering
	\begin{tikzpicture}
		\clip(0.6,0.6) rectangle (4.4,4.4);
		\foreach \a in {(0,0),(0,2),(2,2),(1,1),(2,0),(4,0),(0,4),(2,4),(4,2),%
						(4,4),(3,3),(3,1),(1,3)}
			{
			\begin{scope}[shift={(\a)}]
				\node[fvertex] (a) at (0,0) {};
				\node[fvertex] (b) at (1,0) {};
				\node[fvertex] (c) at (1,1) {};
				\node[fvertex] (d) at (0,1) {};
				\draw[edge] (a) to (b) (b) to (c) (c) to (d) (d) to (a)
							(a) to (c) (b) to (d);
			\end{scope}
			}
	\end{tikzpicture}
	\qquad
	\begin{tikzpicture}
		\clip(0.6,0.6) rectangle (4.4,4.4);
		\pgfmathsetmacro\ang{10};
		\pgfmathsetmacro\up{sqrt(2*(1-cos(180-2*\ang)))};

		\foreach \a in {(0,0),(0,2),(2,2),(2,0),(4,0),(0,4),(2,4),(4,2),
						(4,4)}
			{
		    \node[coordinate] (aux) at \a {};
 			\coordinate (t) at ($0.5*\up*(aux)$);
			\begin{scope}[shift={(t)}]
				\begin{scope}[rotate=-\ang]
					\node[fvertex] (a) at (0,0) {};
					\node[fvertex] (b) at (1,0) {};
					\node[fvertex] (c) at (1,1) {};
					\node[fvertex] (d) at (0,1) {};
					\draw[edge] (a) to (b) (b) to (c) (c) to (d) (d) to (a)
								(a) to (c) (b) to (d);
				\end{scope}
			\end{scope}
			}
		\foreach \a in {(1,1),(3,3),(3,1),(1,3)}
			{
			\coordinate (o) at (0,0);
 			\coordinate[rotate around=-\ang:(o)] (baseShift) at (1,1);
		    \node[coordinate] (aux) at \a {};
 			\coordinate (t) at ($0.5*\up*($(aux)-(1,1)$)+(baseShift)$);	
			\begin{scope}[shift={(t)}]
				\begin{scope}[rotate=\ang]
					\node[fvertex] (a) at (0,0) {};
					\node[fvertex] (b) at (1,0) {};
					\node[fvertex] (c) at (1,1) {};
					\node[fvertex] (d) at (0,1) {};
					\draw[edge] (a) to (b) (b) to (c) (c) to (d) (d) to (a)
								(a) to (c) (b) to (d);
				\end{scope}
			\end{scope}
			}
	\end{tikzpicture}
	\caption{A flexible infinite framework with periodic symmetry.}
	\label{fig:chessboard}
\end{figure}

A \emph{realization} of a (simple) graph $G=(V_G,E_G)$ in the plane is a map $\rho:V_G \rightarrow \mathbb{R}^2$ associating to each vertex a point in the plane such that adjacent vertices are mapped to distinct points, i.e.~$\rho(u) \neq \rho(v)$ for each edge $uv \in E_G$.
A graph together with one of its realizations is called a \emph{planar framework},
or simply a \emph{framework} for short.
If a framework can be continuously deformed in the plane keeping the distances between adjacent vertices,
then we say that such framework is \emph{flexible}; otherwise, the framework is called \emph{rigid} (see \Cref{def:flexible}).
An example of a flexible framework with a countably infinite set of vertices can be seen in \Cref{fig:chessboard}.

Rigidity and flexibility of finite frameworks are well-studied topics with a vast amount of literature between them.
For example, the finite graphs that have a rigid generic realization\footnote{A realization in which the coordinates of the points are algebraically independent.} (and equivalently, only have rigid generic realizations) are combinatorially characterized by the so-called \emph{Laman's theorem},
first proven by Pollaczek-Geiringer in 1924 \cite{Geiringer1927} and independently proven by Laman over forty years later \cite{Laman1970}.
However, not much is known on rigidity for countably infinite graphs, i.e.~graph with a countably infinite amount of vertices.
The combinatorial characterization of Pollaczek-Geiringer and Laman
was extended to countably infinite frameworks by Kitson and Power in~\cite{Kitson2018},
where it was shown that a countably infinite graph will have only rigid generic realizations\footnote{A realization in which every restriction of the realization to a finite vertex set is generic.} if and only if it contains an increasing sequence of finite graphs with rigid realizations.
Countably infinite graphs can, however, be constructed with both rigid and flexible generic realizations; for example, see \cite[Figure~4.2]{dewarthesis}.
This raises the question: which rigidity and flexibility results for finite graphs can be extended to infinite graphs, and which cannot?
This question has been looked into in recent years, with particular interest regarding periodic frameworks \cite{Borcea2010, Malestein2013}.

In the case of finite graphs, it was proven by Grasegger, Legersk\'y and Schicho in~\cite{flexibleLabelings}
that the existence of a flexible (possibly not generic) realization of a connected finite graph is equivalent to the existence of a particular type of red-blue coloring of the edges of the graph called a \emph{NAC-coloring}  (see \Cref{def:nac}).
We shall prove that this characterization can indeed be extended to countably infinite graphs.

\begin{theorem}\label{mainthm}
	A countably infinite connected graph admits a flexible realization if and only if it has a NAC-coloring.
\end{theorem}

The proof of the main theorem for finite graphs (see \cite[Theorem~3.1]{flexibleLabelings}) relies on two different techniques.
The easier implication involves constructing a flexible framework from a NAC-coloring,
which can be extended to countably infinite graphs with relative ease (\Cref{proposition:construction}).
However, the algebraic methods required to obtain a NAC-coloring from a finite flexible framework break down for infinite frameworks,
although the technique was adapted to infinite frameworks with periodic symmetry by Dewar in~\cite{Dewar}.
We will bypass this obstacle by instead considering an ascending tower of finite subgraphs
that covers the whole graph and then recursively applying the known result about NAC-colorings for finite graphs.
In particular, once a set of NAC-colorings with some specific properties is found for each subgraph of the tower,
there is an infinite chain of NAC-colorings extending each other by K\"onig's lemma (\Cref{lemma:koenig}).
This infinite chain will give us a NAC-coloring of the infinite graph.
\Cref{mainthm} was originally presented at \cite{GLSMEGA},
however the proof technique was vastly different and
utilized much more heavy machinery from algebraic geometry; see \Cref{rem:MEGAapproach}.

There is a particular type of infinite framework that we shall wish to apply \Cref{mainthm} and its methods to: the \emph{Penrose frameworks}, i.e.~frameworks that are the 1-skeletons of Penrose rhombus tilings (see \Cref{sec:penrose} for more details).
The eponymous tiling was first formulated by Penrose in the 1970s \cite{penrose74} and is one of the first examples of a \emph{quasi-crystal},
an infinite structure that has order but is not periodic.
Ideally, we would like to be able to determine which rhombi we should brace to obtain rigidity solely from the \emph{ribbon graph} of the Penrose framework (\Cref{def:ribbongraph}).
A major hurdle when working with Penrose frameworks is that they are neither periodic nor generic, so most previous tools cannot be used with them.

To determine which braced Penrose frameworks are rigid, we shall study the much larger class of \emph{P-frameworks}  (\Cref{def:pframework}).
First defined in \cite{Grasegger2020} as the generalization of 1-skeletons of parallelogram tilings,
it was shown by Grasegger and Legersk\'y for a finite P-framework $(G,\rho)$ with some parallelograms braced,
the following properties are equivalent: (i)  $(G,\rho)$ is rigid, (ii) $G$ has no \emph{cartesian NAC-coloring} (\Cref{def:cartesian}), and (iii) the \emph{bracing graph} (\Cref{def:ribbongraph}) of $(G, \rho)$ is connected.
We shall extend this result to countably infinite P-frameworks, allowing us to give a characterization for rigid braced Penrose frameworks. 

\begin{theorem}
	\label{thm:bracedPframeworks}
	Let $(G,\rho)$ be a countably infinite braced P-framework.
	Then the following are equivalent:
	\begin{enumerate}[(i)]
		\item $(G,\rho)$ is rigid. \label{it:rigid}
		\item $G$ has no cartesian NAC-coloring. \label{it:noNAC}
		\item The bracing graph of $G$ is connected. \label{it:connected}
	\end{enumerate}
\end{theorem}

As some Penrose tilings have 5-fold rotational symmetry,
we would also like to be able to determine whether there exists a flex of a Penrose framework that maintains the 5-fold rotational symmetry.
To do this, we shall first apply the methods used in the proof of \Cref{mainthm} to infinite graphs with $k$-fold rotational symmetry to give a combinatorial characterization for the existence of $k$-fold rotationally symmetric flexible realizations.
This will use methods employed by Dewar, Grasegger and Legersk\'y,
who gave a combinatorial characterization of finite rotationally symmetric graphs with flexible symmetric realizations~\cite{RotSymmetry}.
We shall then combine this result with the methods used to prove \Cref{thm:bracedPframeworks} to obtain a combinatorial characterization of P-frameworks with rotationally symmetric flexible motions.
With this we will be able to prove the following result.

\begin{corollary}\label{cor:SymmPenrose}
	Any Penrose framework with 5-fold rotational symmetry has a non-trivial flex (\Cref{def:flexible})
	that preserves 5-fold rotational symmetry.
\end{corollary}

\begin{figure}
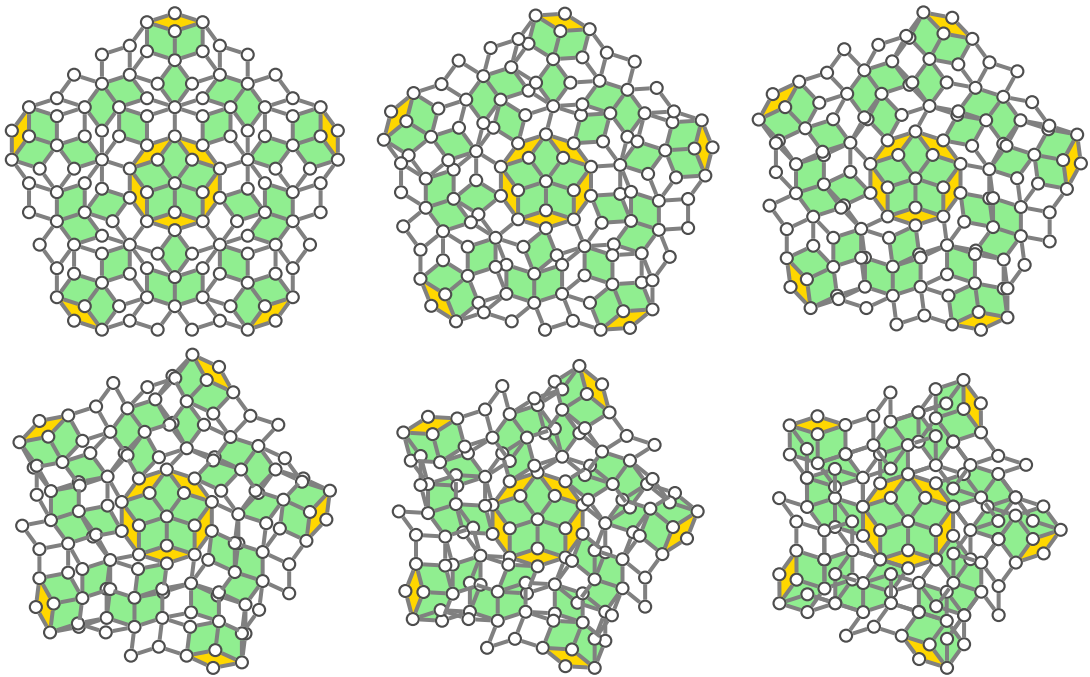

	\centering


	\caption{A 5-fold rotationally symmetric deformation of a Penrose framework. The filled rhombi preserve their shapes during the flex.}
	\label{fig:flexingPenrose}
\end{figure}

See \Cref{fig:flexingPenrose} for an example of one of the 5-fold rotationally symmetric Penrose frameworks
with a 5-fold rotationally symmetric 
flex\footnote{An animation of the flex together with the code producing it using package \FlexRiLoG{} \cite{flexrilog}
can be found at \url{https://jan.legersky.cz/PenroseFrameworks}}.

The paper shall be organized as follows.
\Cref{preliminaries} will cover all the required background for rigid and flexible frameworks, NAC-colorings, and all of the algebraic and combinatorial tools we shall require throughout.
In \Cref{general}, we will prove our main result of \Cref{mainthm}.
\Cref{sec:braced} will introduce P-frameworks and will provide a proof for \Cref{thm:bracedPframeworks}.
\Cref{sec:penrose} will formally define Penrose tilings and frameworks, and will also contain the corollaries we can obtain by applying \Cref{thm:bracedPframeworks}.
\Cref{sec:rot} will introduce rotational symmetry for frameworks in general, and will also extend the main result of \cite{RotSymmetry} to frameworks with a countably infinite amount of vertices,
and \Cref{sec:rot+braced} will combine the results of \Cref{sec:braced,sec:rot} to give a characterization of $k$-fold rotationally symmetric braced P-frameworks with flexes that preserve the symmetry.
Finally, in \Cref{sec:dixon} we shall describe all the possible realizations of the countably infinite bipartite framework that are both injective and flexible.
This is notable as our construction detailed in \Cref{proposition:construction} cannot form injective realizations of this graph.

\section{Preliminaries}
\label{preliminaries}

We recall the notions of rigidity and flexibility, 
which are usually given for finite graphs, but extend naturally to graphs with countably many vertices and edges.
We denote by $(\lambda_e)_{e\in E_G}$, where $\lambda_{uv}=\|\rho(u) - \rho(v) \|$, the edge lengths induced by a framework $(G,\rho)$.
For $\KK \in \{\RR, \CC\}$,
we will denote the set of maps $\rho: V_G \rightarrow \KK^2$ by $(\KK^2)^{V_G}$.
When $G$ is finite, we will make no distinction between $(\KK^2)^{V_G}$ and $\KK^{2|V_G|}$.

\begin{definition}\label{def:flexible}
	A \emph{flex} of the framework $(G,\rho)$ is a continuous path $t \mapsto \rho_t$, $t \in [0,\delta)$ for some $\delta>0$,
	in the space of realizations of~$G$\footnote{The topology we will gift the set $(\RR^2)^{V_G}$ will always be the product topology, 
	hence a path $t \mapsto \rho_t$ will be continuous if and only if each path $t \mapsto \rho_t(v)$ is continuous.}
	such that $\rho_0= \rho$ and each $(G,\rho_t)$
	induces the same edge lengths as~$(G,\rho)$.
	The flex is called \emph{trivial} if for all $t \in [0,\delta)$ the realization $\rho_t$
	is congruent to $\rho$, namely, there exists a Euclidean
	isometry~$M_t$ of $\RR^2$ such that $M_t \rho_t(v) = \rho(v)$ for all $v \in V_G$.
	
	We define a framework to be \emph{flexible}
	if there is a non-trivial flex in $\RR^2$.
	Otherwise, it is called \emph{rigid}.
\end{definition}

Due to the nature of our methods, we will want to transform the information about the existence of a flex
into the existence of a positive-dimensional algebraic set consisting of realizations.
To do so, we will need to remove all but finitely many congruent copies of every realization.
The usual method of doing this is to add extra algebraic constraints that stop any rigid body motions of the framework, such as requiring that a single edge is fixed in place.
In order to deal with various special cases like symmetry, we introduce the following more general notion.

\begin{definition}
	Let $H$ be a finite graph.
	Let $X$ be an algebraic set in $(\CC^{2})^{V_H}$ given by polynomials with real coefficients
	such that for every $\rho\in (\RR^{2})^{V_H} \cap X$ 
	there are only finitely many realizations $\rho'\in (\RR^{2})^{V_H} \cap X$ congruent to $\rho$.
	We call any such set $X$ a \emph{restricting set}.
	We use variables $(x_v,y_v)$ for the coordinates of vertex $v\in V_H$.
\end{definition}
\begin{example}
	\label{ex:fixedEdge}
	Let $(H,\rho)$ be a finite framework.
	Let $\bar{u}\bar{v}$ be an edge of~$H$ and $(\lambda_e)_{e\in E_H}$ be the edge lengths induced by $\rho$.
	The zero set of the equations
	\begin{equation*}
		\begin{aligned}
			x_{\bar{u}}=0\,, \quad	y_{\bar{u}}=0\,, \quad
			x_{\bar{v}}=\lambda_{\bar{u}\bar{v}}\,, \quad y_{\bar{v}}=0\,.
		\end{aligned}
	\end{equation*}
	is a restricting set.
	Since these four equations impose that
	the vertices of the edge $\bar{u}\bar{v}$ have a fixed position,
	there is only one realization congruent to a given one, namely, its reflection w.r.t.\ the $x$-axis. 
\end{example}
Later in the paper, we use other restricting sets to deal with rotationally symmetric frameworks and frameworks with parallelograms.

\begin{definition}
	\label{def:mainSystemOfEquations}
	Let $H$ be a finite graph and let $(H,\rho)$ be a flexible framework with induced edge lengths~$(\lambda_e)_{e\in E_H}$.
	Let $X$ be a restricting set.
	We consider the following polynomial system in~$(\CC^2)^{V_H}$ for unknown coordinates $(x_u,y_u)$ where $u\in V_H$:
	\begin{align} \label{eq:mainSystemOfEquations}
		(x_u-x_v)^2+(y_u-y_v)^2&= \lambda_{uv}^2 \quad \text{ for all } uv \in E_H \,.
	\end{align}
	An irreducible complex algebraic curve~$\C$ in the intersection of $X$
	and the zero set of~\eqref{eq:mainSystemOfEquations}
	is called an \emph{algebraic motion of $(H,\rho)$ w.r.t.\ $X$}.
	If the restricting set $X$ is the one obtained as in \Cref{ex:fixedEdge},
	we call it an \emph{algebraic motion of $(H,\rho)$ w.r.t.\ $\bar{u}\bar{v}$}.
\end{definition}

Note that if $(H, \rho)$ has a non-trivial flex, then the system
given by~\eqref{eq:mainSystemOfEquations} and \Cref{ex:fixedEdge} has infinitely
many solutions and so $(H, \rho)$ admits an algebraic motion w.r.t.\ $\bar{u}\bar{v}$.

The techniques employed in~\cite{flexibleLabelings} prove that a finite connected graph has a flexible framework
if and only if it admits a so-called NAC-coloring.
We extend the notion of NAC-colorings to countable infinite graphs in a straightforward manner.

\begin{definition}\label{def:nac}
	Let~$G$ be a graph (finite or countable). 
	A coloring of edges $\delta\colon E_G \rightarrow \{\text{\blue{}, \red{}}\}$ 
	is called a \emph{NAC-coloring},
	if it is surjective and for every cycle\footnote{Although the graph can be infinite, cycles are always finite.} in~$G$,
	either all edges have the same color, or
	there are at least two edges in each color.
	A NAC-coloring~$\delta$ induces two subgraphs of~$G$:
	\begin{align*}
		G^\delta_\red = (V_G, \{e\in E_G \colon \delta(e)=\red \}) \quad \text{ and } \quad
		G^\delta_\blue = (V_G, \{e\in E_G \colon \delta(e)=\blue \})\,.
	\end{align*}
	The set of NAC-colorings of~$G$ is denoted by~$\nac{G}$.
	See \Cref{fig:NACchessboard} for an example.
\end{definition}

\begin{figure}[ht]
	\centering
	\begin{tikzpicture}
		\clip(0.6,0.6) rectangle (4.4,4.4);
		\foreach \a in {(0,0),(0,2),(2,2),(2,0),(4,0),(0,4),(2,4),(4,2),
						(4,4)}
			{
			\begin{scope}[shift={(\a)}]
				\node[fvertex] (a) at (0,0) {};
				\node[fvertex] (b) at (1,0) {};
				\node[fvertex] (c) at (1,1) {};
				\node[fvertex] (d) at (0,1) {};
				\draw[redge] (a) to (b) (b) to (c) (c) to (d) (d) to (a)
							(a) to (c) (b) to (d);
			\end{scope}
			}
		\foreach \a in {(1,1),(3,3),(3,1),(1,3)}
			{
			\begin{scope}[shift={(\a)}]
				\node[fvertex] (a) at (0,0) {};
				\node[fvertex] (b) at (1,0) {};
				\node[fvertex] (c) at (1,1) {};
				\node[fvertex] (d) at (0,1) {};
				\draw[bedge] (a) to (b) (b) to (c) (c) to (d) (d) to (a)
							(a) to (c) (b) to (d);
			\end{scope}
			}
	\end{tikzpicture}
	\caption{A NAC-coloring of the graph in \Cref{fig:chessboard}}
	\label{fig:NACchessboard}
\end{figure}

As was shown in \cite{flexibleLabelings}, all NAC-colorings can be constructed starting from valuations on the function field (see \cite{Stichtenoth08}) of an algebraic motion $\C$. 
The following functions are used:
	\begin{align*}
		W_{u,v}(\rho) := (x_v - x_u) + i (y_v - y_u), \qquad Z_{u,v}(\rho) := (x_v - x_u) - i (y_v - y_u),
	\end{align*}
where $uv$ is an edge of the graph and $(x_u,y_u)=\rho(u)$ for $\rho\in\C$.
This allows the assignment of some NAC-colorings to an algebraic motion (see also \cite[Definition~2.9]{movableGraphs}).
\begin{definition}
	\label{def:active}
	Let $H$ be a finite graph and let $(H,\rho)$ be a flexible framework.
	Let~$\C$ be an algebraic motion of~$(H,\rho)$ w.r.t.\ a restricting set $X$.
	A NAC-coloring $\delta$ of $H$ is called \emph{active w.r.t.\ $\C$} if there exists a valuation $\nu$ on the function field of~$\C$
	and $\beta \in \QQ$ such that $\delta(uv)=\red \iff \nu(W_{u,v})>\beta$ for $uv\in E_H$.
\end{definition}

We finish the section with an important technical lemma we will require throughout the paper.

\begin{lemma}\label{lem:active}
	Let $(H,\rho)$ be a finite framework and $u_1v_1,u_2v_2 \in E_H$.
	Let $X$ be a restricting set.
	Let $\alpha$ be a non-trivial flex of $(H,\rho)$ such that for every realization in the flex
	there is a congruent realization in $X$.
	If the angle between the edges $u_1v_1$ and $u_2v_2$ is non-constant along the flex $\alpha$,
	then there exists an algebraic motion $\C$ of $(H,\rho)$ w.r.t. $X$
	and an active NAC-coloring $\delta$ w.r.t. $\C$ such that $\delta(u_1v_1)\neq \delta(u_2v_2)$.
\end{lemma}
\begin{proof}
	Let $A\subseteq X$ be the algebraic set of frameworks with the same induced edge lengths as $\rho$.
	For every realization in the flex~$\alpha$, there is a congruent one contained in $A$ by assumption.
	Define $f:A \rightarrow \CC$ to be the map
	\begin{align*}
		f(\rho') = (\rho'(u_1)-\rho'(v_1))\cdot(\rho'(u_2)-\rho'(v_2)).
	\end{align*}
	As algebraic sets only have finitely many irreducible components and
	there are infinitely many values of the angle between $u_1v_1$ and $u_2v_2$ along the flex $\alpha$,
	there exists an irreducible component $A_0$ of $A$
	which contains realizations $\rho_1,\rho_2 \in A_0$ where $f(\rho_1) \neq f(\rho_2)$.
	By \cite[Lemma pg.~56]{mumford},
	there exists an algebraic curve $\mathcal{C}$ in $A_0$ that contains $\rho_1$ and $\rho_2$.
	
	It is immediate that $W_{u,v}Z_{u,v} =\lambda_{uv}^2$ for all $uv \in E_H$.
	Since
	\begin{equation}
		2f = W_{u_1,v_1}Z_{u_2,v_2} + W_{u_2,v_2}Z_{u_1,v_1} = 
			 W_{u_1,v_1}\frac{\lambda_{u_2v_2}^2}{W_{u_2,v_2}} + W_{u_2,v_2}\frac{\lambda_{u_1v_1}^2}{W_{u_1,v_1}}\,,
	\end{equation}
	the map $W_{u_1,v_1}/W_{u_2,v_2}$ is not constant on $\mathcal{C}$.
	By Chevalley's theorem, there is a valuation~$\nu$ of the complex function field of $\C$
	such that $\nu(W_{u_1,v_1}/W_{u_2,v_2})>0$, hence, $\nu(W_{u_1,v_1})>\nu(W_{u_2,v_2})$.
	We now apply the methods outlined in \cite[Theorem 3.1]{flexibleLabelings}
	(see also \cite[Theorem~2.8]{movableGraphs})
	to construct an active NAC-coloring~$\delta$ w.r.t\ $\C$ according to \Cref{def:active} with $\beta:=\nu(W_{u_2,v_2})$.
	By construction, we have $\delta(u_1v_1) \neq \delta(u_2v_2)$.
\end{proof}

\section{Flexible realizations of infinite graphs}
\label{general}

In this section we shall prove \Cref{mainthm}.
To begin with, we shall prove that every NAC-coloring yields a flexible framework.
The following statement was proved for the case of finite graphs in \cite[Theorem~3.1]{flexibleLabelings}.

\begin{proposition}
\label{proposition:construction}
	Let $G$ be a countably infinite graph.
	If $G$ has a NAC-coloring, then there is a flexible framework $(G,\rho)$.
\end{proposition}
\begin{proof}
	Let~$\delta$ be a NAC-coloring of $G$.
	Let $R_1, R_2, \dots$ be the sets of vertices of connected components of the graph $G^\delta_\red{}$ and
	$B_1, B_2, \dots$ be the sets of vertices of connected components of the graph $G^\delta_\blue{}$.
	Let $r_1, r_2 \dots$ and $b_1, b_2, \dots$ be pairwise distinct vectors in~$\RR^2$
	such that $r_1 = b_1 = (0,0)$.
	For $t \in {[0, 2\pi )}$,
	we define a map $\rho_t\colon V_G\rightarrow \RR^2$ by
	\begin{equation*}
	  \rho_t(v)=
	  \begin{pmatrix}
			\cos t & \sin t \\
			-\sin t & \cos t
		\end{pmatrix}%
		\cdot b_j +r_i \,,
	\end{equation*}
	where $i$ and $j$ are the unique indices such that~$v\in R_i \cap B_j$.
	It can now be checked that $(G, \rho) = (G,\rho_0)$ is a realization (i.e.~$\rho(u) \neq \rho(v)$ for all $uv \in E$) and $t \mapsto \rho_t$ is a flex of~$(G,\rho)$;
	see \cite[Theorem~3.1]{flexibleLabelings} for more details.
\end{proof}

A useful technique when working with countably infinite graphs is to consider increasing sequences of subgraphs.
We define this more explicitly as follows.
\begin{definition}
	Let $G$ be a connected graph. 
	A sequence $(G_n)_{n \in \NN}$ of induced finite connected subgraphs of~$G$ is a \emph{subgraph tower} in~$G$
	if $G_n$ is a proper subgraph of $G_{n+1}$ for all $n \in \NN$, and $\bigcup_{n \in \NN} V_{G_n} = V_G$.\footnote{Our definition differs slightly from that given in \cite{Kitson2018}, as we shall always assume that all subgraphs in a subgraph tower are induced and the tower is vertex spanning.}
\end{definition}

The construction of a NAC-coloring of an infinite graph from a sequence of NAC-colorings of a subgraph tower is
inspired by the proof of De Bruijn-Erdős theorem for countable graphs, see for instance~\cite{NashWilliams1967}.
We use the following version of König's infinity lemma.
\begin{lemma}[\cite{Koenig1927}]
	\label{lemma:koenig}
	Let $(S_n)_{n \in \NN}$ be an infinite sequence of disjoint non-empty finite sets.
	Let $\prec$ be a relation in $\bigcup_{n \in \NN} S_n$ such that for every positive integer $n$
	and $x \in S_{n+1}$, there exists $y\in S_n$ with $y\prec x$. 
	There exists an infinite sequence $(x_n)_{n \in \NN}$ such that $x_n \in S_n$ and $x_n \prec x_{n+1}$
	for all positive integers $n$.
\end{lemma}

\begin{lemma}
	\label{lem:infNACfromSequence}
	Let $G$ be a countably infinite graph with a subgraph tower $(G_n)_{n \in \NN}$.
	Let $e_1,e_2$ be two edges of $G_1$.
	If for all $n\in \NN$, there is a NAC-coloring $\delta_n$ of $G_n$ such that $\delta_n(e_1)\neq\delta_n(e_2)$,
	then $G$ has a NAC-coloring.	
\end{lemma}
\begin{proof}
	Let  
	\begin{equation*}
		S_n = \{\delta\in\nac{G_n} : \delta(e_1)=\blue \neq \delta(e_2)\} 
	\end{equation*}
	for every $n \in \NN$.
	Each $S_n$ is non-empty since either $\delta_n$ is in $S_n$,
	or the NAC-coloring obtained by swapping the colors from $\delta_n$ is in $S_n$. 
	Let $\prec$ be the relation in $\bigcup_{n \in \NN} S_n$ given by restriction to edges of the smaller graph.
	Clearly, a restriction of a NAC-coloring to a subgraph is a NAC-coloring,
	provided that the obtained coloring is surjective,
	which is guaranteed by the condition $\delta(e_1)=\blue \neq \delta(e_2)$.
	Hence, the assumption of \Cref{lemma:koenig} holds.
	Therefore, we get a sequence of NAC-colorings of $(G_n)_{n \in \NN}$ extending each other.
	
	This sequence of NAC-colorings gives a well-defined coloring for~$G$,
	since each edge in~$G$ is eventually an edge of a subgraph in the tower. 
	Such a coloring is a NAC-coloring since every cycle of~$G$ is contained in some $G_j$ in the subgraph tower.
\end{proof}

\begin{proof}[Proof of \Cref{mainthm}]
	If $G$ has a NAC-coloring then there is a flexible realization of $G$ by \Cref{proposition:construction}.
	Now suppose $G$ has a flexible realization $\rho$ with a non-trivial flex $\alpha$.
	Since $\alpha$ is non-trivial, there are edges $e_1,e_2 \in E_G$ whose angle changes along the flex.
	We choose any subgraph tower $(G_n)_{n \in \mathbb{N}}$ such that $e_1,e_2 \in E_{G_n}$ for each $n \in \NN$.
	Clearly, the restriction of the flex $\alpha$ to $V_{G_n}$ is a flex of $G_n$
	such that the angle between $e_1$ and $e_2$ changes.   
	Hence by \Cref{lem:active} with the restricting set given in \Cref{ex:fixedEdge},
	there is a NAC-coloring $\delta_n$ of $G_n$ such that $\delta_n(e_1)\neq\delta_n(e_2)$ for all $n\in \NN$.
	A NAC-coloring of $G$ is now provided by \Cref{lem:infNACfromSequence}.
\end{proof}

\begin{remark}
	\label{rem:MEGAapproach}
	As we mentioned in the introduction,
	\Cref{mainthm} was originally presented at \cite{GLSMEGA} with the following proof technique.
	A chain of algebraic motions restricting to each other was constructed for a subgraph tower.
	With this, an active NAC-coloring of the first motion can be extended to an active NAC-coloring of the second etc.\	and
	this is used to obtain a NAC-coloring of the infinite graph.
\end{remark}

\section{Bracing of infinite P-frameworks}\label{sec:braced}

We will apply the methods used in the previous section to a specific class of frameworks that was first introduced in \cite{Grasegger2020}.
We begin with the following definitions.

\begin{definition}
	Let $G$ be a graph. Consider the relation on the set of edges, where
	two edges are in relation if they are opposite edges of a 4-cycle subgraph of $G$.
	An~equivalence class of the reflexive-transitive closure of the relation is called a \emph{ribbon}.
	A ribbon $r$ is \emph{simple} if the subgraph induced by $r$ does not contain any 4-cycle.
\end{definition}

\begin{definition}
	Let $G$ be a connected graph.
	A realization $\rho:V_G\rightarrow \RR^2$ for $G$ such that $\rho$ 
	is injective and each 4-cycle in $G$ forms a non-degenerate (i.e., not all vertices are collinear) parallelogram in $\rho$ 
	is called a \emph{parallelogram realization}.
\end{definition}

\begin{definition}\label{def:pframework}
	A graph is called \emph{ribbon-cutting graph} if it is connected and every ribbon is an edge cut,
	namely, removal of the edges of the ribbon disconnects the graph.
	If $\rho$ is a parallelogram realization of a ribbon-cutting graph $G$, we call the framework $(G,\rho)$ a \emph{P-framework}.
\end{definition}

\begin{remark}
	Although the parallelograms in a P-framework are required to be non-degenerate,
	we do not impose that a non-trivial flex of the P-framework consists of parallelogram realizations.
\end{remark}

As we can see with the next result, the existence of a parallelogram realization of a ribbon-cutting graph
will imply certain structural properties for the graph. 

\begin{proposition}
	\label{prop:PframeworkProperties}
	Let $(G,\rho)$ be a P-framework. All ribbons of $G$ are simple
	and every two vertices of $G$ are separated by a ribbon.
\end{proposition}
\begin{proof}
	The proof is analogous to the one for the finite case, see \cite[Remark~3.5, Theorem~3.9]{Grasegger2020}.
\end{proof}

\begin{definition}\label{def:ribbongraph}
	A \emph{braced ribbon-cutting graph} is a graph $G=(V_G,E_c\cup E_d)$ where $E_c$ and $E_d$ are two non-empty disjoint sets 
	such that the graph $(V_G,E_c)$ is a ribbon-cutting graph
	and the edges in $E_d$ correspond to diagonals of some 4-cycles of $(V_G,E_c)$.
	These diagonals are also called \emph{braces}.
	If $r$ is a ribbon of $(V_G,E_c)$, then 
	\begin{equation*}
		r \cup \{u_1u_3 \in E_d \colon \exists \text{ 4-cycle } (u_1,u_2,u_3,u_4) \text{ of } (V_G,E_c) \text{ s.t. }  u_1u_2,u_3u_4\in r\}
	\end{equation*}
	 is a \emph{ribbon} of the braced ribbon-cutting graph~$G$.
	
	The framework $(G,\rho)$ is called \emph{braced P-framework} 
	if $G$ is a braced ribbon-cutting graph and $\rho$ is a parallelogram realization for $(V_G,E_c)$.
\end{definition}

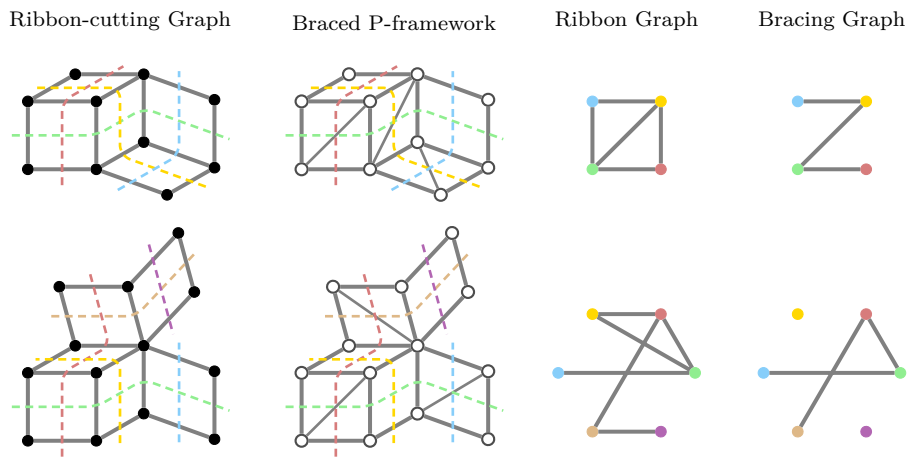
\begin{figure}[ht]
  \centering
  \begin{tikzpicture}[scale=0.9]
		\begin{scope}[yshift=1.9cm]
		  \node[align=center,anchor=south,font=\scriptsize] at (1.35,0) {Ribbon-cutting Graph};
		  \node[align=center,anchor=south,font=\scriptsize] at (5.35,0) {Braced P-framework};
		  \node[align=center,anchor=south,font=\scriptsize] at (8.75,0) {Ribbon Graph};
		  \node[align=center,anchor=south,font=\scriptsize] at (11.75,0) {Bracing Graph};
		\end{scope}

		\begin{scope}
			\begin{scope}
				\node[gvertex] (a) at (0,0) {};
				\node[gvertex] (b) at (1,0) {};
				\node[gvertex] (c) at (1,1) {};
				\node[gvertex] (d) at ($(a)+(c)-(b)$) {};
				\node[gvertex,rotate around=30:(b)] (e) at ($(b)+(0.8,0)$) {};
				\node[gvertex] (f) at ($(e)+(c)-(b)$) {};
				\node[gvertex,rotate around=-20:(e)] (g) at ($(e)+(1.1,0)$) {};
				\node[gvertex] (h) at ($(g)+(f)-(e)$) {};
				\node[gvertex] (j) at ($(b)+(g)-(e)$) {};
				\node[gvertex] (k) at ($(d)+(f)-(c)$) {};
				
				\draw[edge] (a)--(b) (b)--(c) (c)--(d) (d)--(a) (b)--(e) (e)--(f) (f)--(c) (e)--(g) (g)--(h) (h)--(f) (g)--(j) (b)--(j) (d)--(k) (f)--(k);
				\draw[ribbon,col1] ($(a)!0.5!(d)$) -- ($(b)!0.5!(c)$) -- ($(e)!0.5!(f)$) -- ($(g)!0.5!(h)$);
				\draw[ribbon,col2] ($(a)!0.5!(b)$) -- ($(d)!0.5!(c)$) -- ($(f)!0.5!(k)$);
				\draw[ribbon,col3] ($(g)!0.5!(j)$) -- ($(b)!0.5!(e)$) -- ($(c)!0.5!(f)$) -- ($(d)!0.5!(k)$);
				\draw[ribbon,col4] ($(b)!0.5!(j)$) -- ($(e)!0.5!(g)$) -- ($(f)!0.5!(h)$);
			\end{scope}
			
			\begin{scope}[xshift=4cm]
				\node[fvertex] (a) at (0,0) {};
				\node[fvertex] (b) at (1,0) {};
				\node[fvertex] (c) at (1,1) {};
				\node[fvertex] (d) at ($(a)+(c)-(b)$) {};
				\node[fvertex,rotate around=30:(b)] (e) at ($(b)+(0.8,0)$) {};
				\node[fvertex] (f) at ($(e)+(c)-(b)$) {};
				\node[fvertex,rotate around=-20:(e)] (g) at ($(e)+(1.1,0)$) {};
				\node[fvertex] (h) at ($(g)+(f)-(e)$) {};
				\node[fvertex] (j) at ($(b)+(g)-(e)$) {};
				\node[fvertex] (k) at ($(d)+(f)-(c)$) {};
				
				\draw[edge] (a)--(b) (b)--(c) (c)--(d) (d)--(a) (b)--(e) (e)--(f) (f)--(c) (e)--(g) (g)--(h) (h)--(f) (g)--(j) (b)--(j) (d)--(k) (f)--(k);
				\draw[brace] (a)--(c) (b)--(f) (e)--(j);
				\draw[ribbon,col1] ($(a)!0.5!(d)$) -- ($(b)!0.5!(c)$) -- ($(e)!0.5!(f)$) -- ($(g)!0.5!(h)$);
				\draw[ribbon,col2] ($(a)!0.5!(b)$) -- ($(d)!0.5!(c)$) -- ($(f)!0.5!(k)$);
				\draw[ribbon,col3] ($(g)!0.5!(j)$) -- ($(b)!0.5!(e)$) -- ($(c)!0.5!(f)$) -- ($(d)!0.5!(k)$);
				\draw[ribbon,col4] ($(b)!0.5!(j)$) -- ($(e)!0.5!(g)$) -- ($(f)!0.5!(h)$);
			\end{scope}
			
			\begin{scope}[xshift=8.25cm]
				\node[gvertex,col1] (1) at (0,0) {};
				\node[gvertex,col2] (2) at (1,0) {};
				\node[gvertex,col3] (3) at (1,1) {};
				\node[gvertex,col4] (4) at (0,1) {};
				\draw[edge] (1)--(2) (2)--(3) (3)--(4) (4)--(1) (1)--(3);
			\end{scope}
			
			\begin{scope}[xshift=11.25cm]
				\node[gvertex,col1] (1) at (0,0) {};
				\node[gvertex,col2] (2) at (1,0) {};
				\node[gvertex,col3] (3) at (1,1) {};
				\node[gvertex,col4] (4) at (0,1) {};
				\draw[edge] (1)--(2) (3)--(4) (1)--(3);
			\end{scope}
		\end{scope}
  
		\begin{scope}[yshift=-4cm]
			\begin{scope}
				\node[gvertex] (a) at (0,0) {};
				\node[gvertex] (b) at (1,0) {};
				\node[gvertex] (c) at (1,1) {};
				\node[gvertex] (d) at ($(a)+(c)-(b)$) {};
				\node[gvertex,rotate around=30:(b)] (e) at ($(b)+(0.8,0)$) {};
				\node[gvertex] (f) at ($(e)+(c)-(b)$) {};
				\node[gvertex,rotate around=-20:(e)] (g) at ($(e)+(1.1,0)$) {};
				\node[gvertex] (h) at ($(g)+(f)-(e)$) {};
				\node[gvertex] (k) at ($(d)+(f)-(c)$) {};
				\node[gvertex,rotate around=105:(k)] (l) at ($(k)+(0.9,0)$) {};
				\node[gvertex] (m) at ($(f)+(l)-(k)$) {};
				\node[gvertex,rotate around=25:(k)] (n) at ($(m)+(1.3,0)$) {};
				\node[gvertex] (o) at ($(f)+(n)-(m)$) {};
				
				\draw[edge] (a)--(b) (b)--(c) (c)--(d) (d)--(a) (b)--(e) (e)--(f) (f)--(c) (e)--(g) (g)--(h) (h)--(f) (d)--(k) (f)--(k) (k)--(l) (l)--(m) (m)--(f) (m)--(n) (n)--(o) (o)--(f);
				\draw[ribbon,col1] ($(a)!0.5!(d)$) -- ($(b)!0.5!(c)$) -- ($(e)!0.5!(f)$) -- ($(g)!0.5!(h)$);
				\draw[ribbon,col2] ($(a)!0.5!(b)$) -- ($(d)!0.5!(c)$) -- ($(f)!0.5!(k)$) -- ($(l)!0.5!(m)$);
				\draw[ribbon,col3] ($(b)!0.5!(e)$) -- ($(c)!0.5!(f)$) -- ($(d)!0.5!(k)$);
				\draw[ribbon,col4] ($(e)!0.5!(g)$) -- ($(f)!0.5!(h)$);
				\draw[ribbon,col5] ($(k)!0.5!(l)$) -- ($(f)!0.5!(m)$) -- ($(n)!0.5!(o)$);
				\draw[ribbon,col6] ($(m)!0.5!(n)$) -- ($(f)!0.5!(o)$);
			\end{scope}
			
			\begin{scope}[xshift=4cm]
				\node[fvertex] (a) at (0,0) {};
				\node[fvertex] (b) at (1,0) {};
				\node[fvertex] (c) at (1,1) {};
				\node[fvertex] (d) at ($(a)+(c)-(b)$) {};
				\node[fvertex,rotate around=30:(b)] (e) at ($(b)+(0.8,0)$) {};
				\node[fvertex] (f) at ($(e)+(c)-(b)$) {};
				\node[fvertex,rotate around=-20:(e)] (g) at ($(e)+(1.1,0)$) {};
				\node[fvertex] (h) at ($(g)+(f)-(e)$) {};
				\node[fvertex] (k) at ($(d)+(f)-(c)$) {};
				\node[fvertex,rotate around=105:(k)] (l) at ($(k)+(0.9,0)$) {};
				\node[fvertex] (m) at ($(f)+(l)-(k)$) {};
				\node[fvertex,rotate around=25:(k)] (n) at ($(m)+(1.3,0)$) {};
				\node[fvertex] (o) at ($(f)+(n)-(m)$) {};
				
				\draw[edge] (a)--(b) (b)--(c) (c)--(d) (d)--(a) (b)--(e) (e)--(f) (f)--(c) (e)--(g) (g)--(h) (h)--(f) (d)--(k) (f)--(k) (k)--(l) (l)--(m) (m)--(f) (m)--(n) (n)--(o) (o)--(f);
				\draw[brace] (a)--(c) (e)--(h) (f)--(l);
				\draw[ribbon,col1] ($(a)!0.5!(d)$) -- ($(b)!0.5!(c)$) -- ($(e)!0.5!(f)$) -- ($(g)!0.5!(h)$);
				\draw[ribbon,col2] ($(a)!0.5!(b)$) -- ($(d)!0.5!(c)$) -- ($(f)!0.5!(k)$) -- ($(l)!0.5!(m)$);
				\draw[ribbon,col3] ($(b)!0.5!(e)$) -- ($(c)!0.5!(f)$) -- ($(d)!0.5!(k)$);
				\draw[ribbon,col4] ($(e)!0.5!(g)$) -- ($(f)!0.5!(h)$);
				\draw[ribbon,col5] ($(k)!0.5!(l)$) -- ($(f)!0.5!(m)$) -- ($(n)!0.5!(o)$);
				\draw[ribbon,col6] ($(m)!0.5!(n)$) -- ($(f)!0.5!(o)$);
			\end{scope}
			
			\begin{scope}[xshift=8.75cm,yshift=1cm]
				\coordinate (o) at (0,0);
				\node[gvertex,col1] (1) at (1,0) {};
				\node[gvertex,col2,rotate around=60:(o)] (2) at (1,0) {};
				\node[gvertex,col3,rotate around=120:(o)] (3) at (1,0) {};
				\node[gvertex,col4,rotate around=180:(o)] (4) at (1,0) {};
				\node[gvertex,col5,rotate around=240:(o)] (5) at (1,0) {};
				\node[gvertex,col6,rotate around=300:(o)] (6) at (1,0) {};
				\draw[edge] (1)--(2) (1)--(3) (1)--(4) (2)--(3) (2)--(5) (5)--(6);
			\end{scope}
			
			\begin{scope}[xshift=11.75cm,yshift=1cm]
				\coordinate (o) at (0,0);
				\node[gvertex,col1] (1) at (1,0) {};
				\node[gvertex,col2,rotate around=60:(o)] (2) at (1,0) {};
				\node[gvertex,col3,rotate around=120:(o)] (3) at (1,0) {};
				\node[gvertex,col4,rotate around=180:(o)] (4) at (1,0) {};
				\node[gvertex,col5,rotate around=240:(o)] (5) at (1,0) {};
				\node[gvertex,col6,rotate around=300:(o)] (6) at (1,0) {};
				\draw[edge] (1)--(2) (1)--(4) (2)--(5);
			\end{scope}
		\end{scope}
  \end{tikzpicture}
  \caption{Two ribbon-cutting graphs with an example of a braced P-framework as well as the corresponding ribbon graph and bracing graph.
  The vertices in the ribbon graph and the bracing graph are colored in correspondence with the indicated ribbons.}
  \label{fig:ribbonbracinggraph}
\end{figure}

\begin{definition}
	Let $G$ be a braced ribbon-cutting graph. The \emph{ribbon graph} $\Gamma$ of $G$
	is the graph with the set of vertices being the set of ribbons
	of $G$ and two ribbons $r_1,r_2$ are adjacent if and only if 
	there is a 4-cycle $(u_1, u_2, u_3, u_4)$ in the underlying unbraced graph of $G$ such that 
	$u_1 u_2, u_3 u_4 \in r_1$ and $u_1 u_4, u_2 u_3 \in r_2$.
	The subgraph $(V_\Gamma, E_b)$ of $\Gamma$, where 
	\begin{equation*}
		E_b = \{r_1r_2\in E_\Gamma \colon r_1 \cap r_2 \text{ is a non-empty subset of braces of } G\}\,,
	\end{equation*}
	 is called the \emph{bracing (sub)graph}; see \Cref{fig:ribbonbracinggraph} for some examples.
\end{definition}

Similar to the last section, we wish to link the existence of flexible motions of P-frameworks to the existence of NAC-colorings with certain properties.
The result about P-frameworks is even stronger, namely, we can decide about flexibility of a given parallelogram realization, not just its existence,
depending on the existence of this special NAC-coloring.

\begin{definition}\label{def:cartesian}
	A NAC-coloring $\delta$ of a graph $G$ is called \emph{cartesian}
	if no two distinct vertices are connected by a red and blue path simultaneously.
\end{definition}

Equivalently, a NAC-coloring $\delta$ is cartesian if and only if
for every connected component $R$ of $G_\red^\delta$ and $B$ of $G_\blue^\delta$,
the intersection of the vertex sets of $R$ and $B$ contains at most one vertex.
Notice that if the construction of a flex in \Cref{proposition:construction}
is applied to a cartesian NAC-coloring, then the realization is injective.

\begin{lemma}
	\label{lem:cartesianIffRibbonsMonochromatic}
	Let $(G,\rho)$ be a braced P-framework.
	A NAC-coloring of $G$ is cartesian if and only if each ribbon of $G$ is monochromatic. 
\end{lemma}
\begin{proof}
	The proof is analogous to the finite case, see \cite[Theorem~3.9, Lemma~4.1]{Grasegger2020}
\end{proof}

We are now ready for the proof of \Cref{thm:bracedPframeworks}.

\begin{proof}[Proof of \Cref{thm:bracedPframeworks}]
	$\neg \ref{it:noNAC} \Rightarrow \neg \ref{it:rigid}$:
	If $G$ has a cartesian NAC-coloring, then $(G,\rho)$ is flexible
	since the constructive method given by \cite[Theorem 4.4]{Grasegger2020}
	can be extended to countably infinite P-frameworks; see also the construction of flex $\rho_t$ in the proof of \Cref{lem:CnCartNACconstruction}.
	It exploits the construction described in the proof of \Cref{proposition:construction}
	with carefully chosen vectors $r_i$ and $b_j$ so that the obtained flex starts at $\rho$. 
	
	$\neg \ref{it:rigid} \Rightarrow \neg \ref{it:noNAC}$:
	Let $(G,\rho)$ be a flexible infinite P-framework.
	Let $\alpha \colon {[0,1)} \rightarrow (\RR^2)^{V_G}$ be the corresponding non-trivial flex.
	Since the flex is non-trivial, there exist edges $e_1$ and~$e_2$
	such that the angle between these two edges changes along the flex.
	Let $(G_n)_{n \in \NN}$ be a subgraph tower of $G$ such that
	$e_1,e_2 \in E_{G_n}$\footnote{It is not clear whether we can choose each $G_n$ to be a ribbon-cutting graph.
	However, our proof does not require this is true.}.
	For any $n \in \NN$, projecting the flex~$\alpha$ to $G_n$ 
	gives a flex $\alpha_n$ of~$G_{n}$.
	We would like to use \Cref{lem:active}
	with the restricting set $X$ given by the equations from \Cref{ex:fixedEdge} and the following
	\begin{equation} \label{eq:parallelograms}
		\begin{aligned}
			x_{u_2} - x_{u_1} &= x_{u_3} - x_{u_4} \,, & x_{u_4} - x_{u_1} &= x_{u_3} - x_{u_2} \,, \\
			y_{u_2} - y_{u_1} &= y_{u_3} - y_{u_4} \,, & y_{u_4} - y_{u_1} &= y_{u_3} - y_{u_2} \,.
		\end{aligned}
	\end{equation}
	for each 4-cycle $(u_1,u_2,u_3,u_4)$ in $G_n$; this forces all 4-cycles to be realized as parallelograms.
	It is possible that some realizations in the flex $\alpha_n$ will not correspond to a congruent copy in $X$,
	since a parallelogram can fold to an antiparallelogram.
	Nevertheless, if a 4-cycle is at some point in the flex realized as an antiparallelogram
	instead of a parallelogram (or with two opposite vertices coinciding in the case of rhombi),
	then, by continuity of the flex, there is a minimum value of parameter $t$ such that all vertices of the 4-cycle are collinear.
	By the assumption that the parallelograms in a P-framework are non-degenerate, this value $t$ is positive.
	Let $\varepsilon_n$ be the minimum of such values over all 4-cycles in~$G_n$ having this behavior;
	since $G_n$ is finite, there are only finitely many 4-cycles, hence $\varepsilon_n>0$.
	Hence, the realizations in the restriction of the flex to $[0,\varepsilon_n)$ have a congruent realization in $X$.
	By \Cref{lem:active}, for every $n \in \NN$, there is a NAC-coloring $\delta_n$ of~$G_n$ 
	such that $\delta_n(e_1)=\blue$ and $\delta_n(e_2)=\red$.
	Moreover, each ribbon in the constructed active NAC-coloring is monochromatic
	since $W_{u_1,u_2} = W_{u_3,u_4}$ for a 4-cycle $(u_1,u_2,u_3,u_4)$
	by~\eqref{eq:parallelograms}.
	Now let  
	\begin{equation*}
		S_n = \{\delta\in\nac{G_n} : \text{ all ribbons are monochromatic in $\delta$ and }
			\delta(e_1)=\blue \neq \delta(e_2)\}\, 
	\end{equation*}
	for every $n \in \NN$.
	By the discussion above, $\delta_n\in S_n$, namely, $S_n$ is non-empty for every $n \in \NN$.
	Let $\prec$ be the relation in $\bigcup_{n \in \NN} S_n$ given by restriction to edges of the smaller graph.
	Clearly, a restriction of a NAC-coloring to a subgraph is a NAC-coloring,
	provided that the obtained coloring is surjective,
	which is guaranteed by the condition $\delta(e_1)=\blue \neq \delta(e_2)$.
	Since a ribbon in a subgraph is a subset of a ribbon of the supergraph,
	the property that ribbons are monochromatic holds for restrictions as well.
	Hence, the assumption of \Cref{lemma:koenig} holds.
	Therefore, we get a sequence of NAC-colorings of $(G_n)_{n \in \NN}$ extending each other.
	
	By the same reason as in \Cref{lem:infNACfromSequence},
	this sequence gives a NAC-coloring for~$G$.
	Moreover, each ribbon of the NAC-coloring is monochromatic,
	otherwise two edges of the ribbon with different colors together 
	with all edges certifying that they are in the same ribbon are eventually contained in some $G_n$ in the same ribbon,
	which contradicts the assumption imposed in the definition of $S_n$.
	Hence, the NAC-coloring is cartesian by \Cref{lem:cartesianIffRibbonsMonochromatic}.
	
	\ref{it:noNAC} $\Leftrightarrow$ \ref{it:connected}:
	The proof of this for finite P-frameworks (see \cite[Theorem 4.5]{Grasegger2020}) extends to countably infinite P-frameworks.
\end{proof}

In \cite{power2021parallelogram},
Power proved that a \emph{parallelogram framework} $(G,\rho)$ (a P-framework that is the 1-skeleton of a parallelogram tiling of the plane) has a connected bracing graph if and only if it is \emph{infinitesimally rigid};
i.e.~for every map $\mu: V_G \rightarrow \mathbb{R}^2$ with
\begin{align*}
	(\rho(u) - \rho(v) ) \cdot (\mu(u) - \mu(v)) = 0  \qquad \text{ for all } uv \in E_G,
\end{align*}
there exists a real-valued $2 \times 2$ skew-symmetric matrix $S$ and a point $q \in \mathbb{R}^2$ so that $\mu(v) = S \rho(v) +q$ for all $v \in V_G$.
While infinitesimal rigidity is a sufficient condition for rigidity for finite frameworks (see \cite{Asimow1978TheRO}),
and is also a sufficient condition for rigidity for generic countably infinite frameworks in the plane (see \cite[Theorem 4.1]{Kitson2018}),
this is not the case for non-generic frameworks;
see \Cref{fig:InfRigidFlexible} for an example of an infinitesimally rigid but flexible framework.
Since all parallelogram frameworks are non-generic,
a weaker version of \Cref{thm:bracedPframeworks} for countably infinite parallelogram frameworks does not directly follow
from Power's result. 
However,
we can obtain the following result by combining Power's result with \Cref{thm:bracedPframeworks}.

\begin{corollary}\label{c:parallelogram}
	A parallelogram framework is rigid if and only if it is infinitesimally rigid.
\end{corollary}
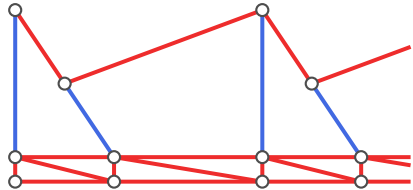
\begin{figure}[ht]
	\centering
	\begin{tikzpicture}[scale = 0.65]
	\draw[redge] (0,0) -- (0,0.5) -- (2,0) -- (2,0.5) -- (5,0) -- (5,0.5) -- (7,0) -- (7,0.5);
	\draw[bedge] (0,0.5) -- (0,3.5);
	\draw[bedge] (5,0.5) -- (5,3.5);
	\draw[redge] (5,3.5) -- (1,2);
	\draw[redge] (6,2) -- (8,2.75);
	\draw[redge] (7,0.5) -- (8,0.3333333333);
	
	\draw[redge] (0,0) -- (2,0) -- (5,0) -- (7,0);
	\draw[redge] (7,0) -- (8,0);
	
	\draw[redge] (0,3.5) -- (1,2);
	\draw[bedge] (1,2) -- (2,0.5);
	\draw[redge] (5,3.5) -- (6,2);
	\draw[bedge] (6,2) -- (7,0.5);
	\draw[redge] (0,0.5) -- (2,0.5) -- (5,0.5) -- (7,0.5);
	\draw[redge] (7,0.5) -- (8,0.5);

	\node[fvertex] (1) at (0,0) {};
	\node[fvertex] (2) at (2,0) {};
	\node[fvertex] (3) at (5,0) {};
	\node[fvertex] (4) at (7,0) {};
	\node[fvertex] (5) at (0,0.5) {};
	\node[fvertex] (6) at (2,0.5) {};
	\node[fvertex] (7) at (5,0.5) {};
	\node[fvertex] (8) at (7,0.5) {};
	\node[fvertex] (9) at (0,3.5) {};
	\node[fvertex] (10) at (5,3.5) {};
	\node[fvertex] (11) at (1,2) {};
	\node[fvertex] (12) at (6,2) {};
	\end{tikzpicture}
	\caption{An infinite framework that is infinitesimally rigid but flexible, as shown in \cite[Section~5]{KASTIS2021125404}.
	Any generic realization of the graph is infinitesimally rigid. The graph has infinitely many NAC-colorings, however the unique cartesian NAC-coloring (up to switching the colors of all of the edges) is the coloring shown. With the construction given in \Cref{proposition:construction} applied to the unique cartesian NAC-coloring, we can obtain an infinite P-framework.}\label{fig:InfRigidFlexible}
\end{figure}

\begin{figure}[ht]
  \centering
  \begin{tikzpicture}[scale=1.2]
    \node[fvertex] (a) at (-0.5,0) {};
    \node[fvertex] (b) at (1.3,0) {};
    \node[fvertex] (c) at (1.3,1) {};
    \node[fvertex] (d) at ($(a)+(c)-(b)$) {};
    \node[fvertex,rotate around=30:(b)] (e) at ($(b)+(0.8,0)$) {};
    \node[fvertex] (f) at ($(e)+(c)-(b)$) {};
    \node[fvertex] (g) at ($(d)+(f)-(c)$) {};
    \node[fvertex] (h) at ($(a)+(g)-(d)$) {};
    
    \draw[edge] (a)--(b) (b)--(c) (c)--(d) (d)--(a) (b)--(e) (e)--(f) (f)--(c) (d)--(g) (f)--(g) (a)--(h) (g)--(h) (e)--(h);
    \draw[ribbon,col1,shorten <= 0pt, shorten >= 0pt] ($(g)!0.5!(h)$) -- ($(a)!0.5!(d)$) -- ($(b)!0.5!(c)$) -- ($(e)!0.5!(f)$) -- cycle;
    \draw[edge] (c)--(d);
    \draw[ribbon,col3,shorten <= 0pt, shorten >= 0pt] ($(e)!0.5!(h)$) -- ($(a)!0.5!(b)$) -- ($(c)!0.5!(d)$) -- ($(f)!0.5!(g)$) -- cycle;
    \draw[ribbon,col4,shorten <= 0pt, shorten >= 0pt] ($(b)!0.5!(e)$) -- ($(c)!0.5!(f)$) -- ($(d)!0.5!(g)$) -- ($(a)!0.5!(h)$) -- cycle;
  \end{tikzpicture}
  \caption{Example of a P-framework with ribbons meeting each other twice.}
  \label{fig:cube}
\end{figure}
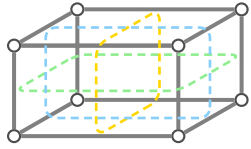

A major obstacle to extending \Cref{c:parallelogram} to P-frameworks is that Power's proof requires
that the union of any two ribbons contains at most one 4-cycle.
While this extra condition holds for parallelogram frameworks, it does not for all P-frameworks;
for example, the union of any two distinct ribbons of 1-skeleton of the cube 
is exactly two disjoint 4-cycles; see \Cref{fig:cube}.

\section{Bracing Penrose frameworks
}\label{sec:penrose}

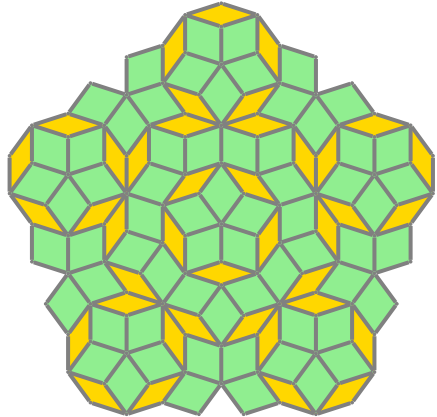
\begin{figure}[ht]
\centering
\begin{tikzpicture}[scale=0.5,rotate=-18]
\fill[thin] (-1.6180, -4.7684e-7) -- (-0.80902, 0.58779) -- (-0.50000, 1.5388) -- (-1.3090, 0.95106) -- cycle;
\fill[thin] (-1.6180, -4.7684e-7) -- (-0.80902, -0.58779) -- (-0.50000, -1.5388) -- (-1.3090, -0.95106) -- cycle;
\fill[thin] (-0.50000, -1.5388) -- (0.30902, -0.95106) -- (1.3090, -0.95106) -- (0.50000, -1.5388) -- cycle;
\fill[thin] (1.3090, -0.95106) -- (1.6180, 1.9073e-6) -- (1.3090, 0.95106) -- (1.0000, 0.00000) -- cycle;
\fill[thin] (1.3090, 0.95106) -- (0.30902, 0.95106) -- (-0.50000, 1.5388) -- (0.50000, 1.5388) -- cycle;
\fill[thin] (-2.6180, -9.5367e-7) -- (-2.3090, 0.95105) -- (-2.6180, 1.9021) -- (-2.9271, 0.95105) -- cycle;
\fill[thin] (-2.6180, -9.5367e-7) -- (-2.9271, 0.95105) -- (-3.7361, 1.5388) -- (-3.4270, 0.58778) -- cycle;
\fill[thin] (-2.6180, -9.5367e-7) -- (-2.9270, -0.95106) -- (-3.7361, -1.5388) -- (-3.4271, -0.58779) -- cycle;
\fill[thin] (-2.6180, -9.5367e-7) -- (-2.3090, -0.95106) -- (-2.6180, -1.9021) -- (-2.9270, -0.95106) -- cycle;
\fill[thin] (-0.80902, -2.4899) -- (-0.49999, -3.4410) -- (0.30902, -4.0287) -- (1.9073e-6, -3.0777) -- cycle;
\fill[thin] (-0.80902, -2.4899) -- (-1.6180, -1.9021) -- (-2.6180, -1.9021) -- (-1.8090, -2.4899) -- cycle;
\fill[thin] (-0.80902, -2.4899) -- (-1.6180, -3.0777) -- (-2.6180, -3.0777) -- (-1.8090, -2.4899) -- cycle;
\fill[thin] (-0.80902, -2.4899) -- (0.19098, -2.4899) -- (1.0000, -3.0777) -- (1.9073e-6, -3.0777) -- cycle;
\fill[thin] (2.1180, -1.5388) -- (3.1180, -1.5388) -- (3.9271, -0.95105) -- (2.9270, -0.95105) -- cycle;
\fill[thin] (2.1180, -1.5388) -- (2.4271, -2.4899) -- (2.1180, -3.4410) -- (1.8090, -2.4899) -- cycle;
\fill[thin] (2.1180, -1.5388) -- (1.3090, -2.1266) -- (1.0000, -3.0777) -- (1.8090, -2.4899) -- cycle;
\fill[thin] (2.1180, -1.5388) -- (2.4271, -0.58778) -- (3.2361, 1.9073e-6) -- (2.9270, -0.95105) -- cycle;
\fill[thin] (2.1180, 1.5388) -- (1.8090, 2.4899) -- (2.1180, 3.4410) -- (2.4270, 2.4899) -- cycle;
\fill[thin] (2.1180, 1.5388) -- (2.9270, 0.95106) -- (3.9270, 0.95106) -- (3.1180, 1.5388) -- cycle;
\fill[thin] (2.1180, 1.5388) -- (2.9270, 0.95106) -- (3.2361, 1.9073e-6) -- (2.4270, 0.58779) -- cycle;
\fill[thin] (2.1180, 1.5388) -- (1.8090, 2.4899) -- (1.0000, 3.0777) -- (1.3090, 2.1266) -- cycle;
\fill[thin] (-0.80902, 2.4899) -- (-1.8090, 2.4899) -- (-2.6180, 3.0777) -- (-1.6180, 3.0777) -- cycle;
\fill[thin] (-0.80902, 2.4899) -- (-1.8090, 2.4899) -- (-2.6180, 1.9021) -- (-1.6180, 1.9021) -- cycle;
\fill[thin] (-0.80902, 2.4899) -- (-1.9073e-6, 3.0777) -- (0.30902, 4.0287) -- (-0.50000, 3.4410) -- cycle;
\fill[thin] (-0.80902, 2.4899) -- (0.19098, 2.4899) -- (1.0000, 3.0777) -- (-1.9073e-6, 3.0777) -- cycle;
\fill[thin] (-3.7361, -1.5388) -- (-4.7361, -1.5388) -- (-5.5451, -0.95106) -- (-4.5451, -0.95106) -- cycle;
\fill[thin] (0.30902, -4.0287) -- (4.7684e-6, -4.9798) -- (-0.80901, -5.5676) -- (-0.50000, -4.6165) -- cycle;
\fill[thin] (3.9271, -0.95105) -- (4.7361, -1.5388) -- (5.0451, -2.4899) -- (4.2361, -1.9021) -- cycle;
\fill[thin] (2.1180, 3.4410) -- (2.9270, 4.0287) -- (3.9270, 4.0287) -- (3.1180, 3.4410) -- cycle;
\fill[thin] (-2.6180, 3.0777) -- (-2.9271, 4.0287) -- (-2.6180, 4.9798) -- (-2.3090, 4.0287) -- cycle;
\fill[thin] (-2.6180, -3.0777) -- (-2.3090, -4.0287) -- (-2.6180, -4.9798) -- (-2.9270, -4.0287) -- cycle;
\fill[thin] (2.1180, -3.4410) -- (2.9271, -4.0287) -- (3.9271, -4.0287) -- (3.1180, -3.4410) -- cycle;
\fill[thin] (3.9270, 0.95106) -- (4.7361, 1.5388) -- (5.0451, 2.4899) -- (4.2361, 1.9021) -- cycle;
\fill[thin] (0.30902, 4.0287) -- (-9.5367e-7, 4.9798) -- (-0.80902, 5.5676) -- (-0.50000, 4.6165) -- cycle;
\fill[thin] (-3.7361, 1.5388) -- (-4.5451, 0.95105) -- (-5.5451, 0.95105) -- (-4.7361, 1.5388) -- cycle;
\fill[thin] (-1.6180, -4.9798) -- (-2.6180, -4.9798) -- (-1.8090, -5.5676) -- (-0.80901, -5.5676) -- cycle;
\fill[thin] (4.2361, -3.0777) -- (5.0451, -2.4899) -- (4.7361, -3.4410) -- (3.9271, -4.0287) -- cycle;
\fill[thin] (4.2361, 3.0777) -- (3.9270, 4.0287) -- (4.7361, 3.4410) -- (5.0451, 2.4899) -- cycle;
\fill[thin] (-1.6180, 4.9798) -- (-0.80902, 5.5676) -- (-1.8090, 5.5676) -- (-2.6180, 4.9798) -- cycle;
\fill[thin] (-5.2361, -1.9073e-6) -- (-5.5451, 0.95105) -- (-5.8541, -9.5367e-7) -- (-5.5451, -0.95106) -- cycle;
\fill[fat] (0.00000, 0.00000) -- (-0.80902, 0.58779) -- (-1.6180, -4.7684e-7) -- (-0.80902, -0.58779) -- cycle;
\fill[fat] (0.00000, 0.00000) -- (0.30902, 0.95106) -- (-0.50000, 1.5388) -- (-0.80902, 0.58779) -- cycle;
\fill[fat] (0.00000, 0.00000) -- (-0.80902, -0.58779) -- (-0.50000, -1.5388) -- (0.30902, -0.95106) -- cycle;
\fill[fat] (0.00000, 0.00000) -- (0.30902, -0.95106) -- (1.3090, -0.95106) -- (1.0000, 0.00000) -- cycle;
\fill[fat] (0.00000, 0.00000) -- (0.30902, 0.95106) -- (1.3090, 0.95106) -- (1.0000, 0.00000) -- cycle;
\fill[fat] (-1.6180, -4.7684e-7) -- (-2.6180, -9.5367e-7) -- (-2.3090, -0.95106) -- (-1.3090, -0.95106) -- cycle;
\fill[fat] (-1.6180, -4.7684e-7) -- (-2.6180, -9.5367e-7) -- (-2.3090, 0.95105) -- (-1.3090, 0.95106) -- cycle;
\fill[fat] (-0.50000, -1.5388) -- (-0.80902, -2.4899) -- (-1.6180, -1.9021) -- (-1.3090, -0.95106) -- cycle;
\fill[fat] (-0.50000, -1.5388) -- (-0.80902, -2.4899) -- (0.19098, -2.4899) -- (0.50000, -1.5388) -- cycle;
\fill[fat] (1.3090, -0.95106) -- (0.50000, -1.5388) -- (1.3090, -2.1266) -- (2.1180, -1.5388) -- cycle;
\fill[fat] (1.3090, -0.95106) -- (2.1180, -1.5388) -- (2.4271, -0.58778) -- (1.6180, 1.9073e-6) -- cycle;
\fill[fat] (1.3090, 0.95106) -- (2.1180, 1.5388) -- (1.3090, 2.1266) -- (0.50000, 1.5388) -- cycle;
\fill[fat] (1.3090, 0.95106) -- (1.6180, 1.9073e-6) -- (2.4270, 0.58779) -- (2.1180, 1.5388) -- cycle;
\fill[fat] (-0.50000, 1.5388) -- (-0.80902, 2.4899) -- (0.19098, 2.4899) -- (0.50000, 1.5388) -- cycle;
\fill[fat] (-0.50000, 1.5388) -- (-0.80902, 2.4899) -- (-1.6180, 1.9021) -- (-1.3090, 0.95106) -- cycle;
\fill[fat] (-1.3090, -0.95106) -- (-2.3090, -0.95106) -- (-2.6180, -1.9021) -- (-1.6180, -1.9021) -- cycle;
\fill[fat] (0.50000, -1.5388) -- (0.19098, -2.4899) -- (1.0000, -3.0777) -- (1.3090, -2.1266) -- cycle;
\fill[fat] (1.6180, 1.9073e-6) -- (2.4271, -0.58778) -- (3.2361, 1.9073e-6) -- (2.4270, 0.58779) -- cycle;
\fill[fat] (0.50000, 1.5388) -- (0.19098, 2.4899) -- (1.0000, 3.0777) -- (1.3090, 2.1266) -- cycle;
\fill[fat] (-1.3090, 0.95106) -- (-2.3090, 0.95105) -- (-2.6180, 1.9021) -- (-1.6180, 1.9021) -- cycle;
\fill[fat] (-2.6180, -9.5367e-7) -- (-3.4270, 0.58778) -- (-4.2361, -9.5367e-7) -- (-3.4271, -0.58779) -- cycle;
\fill[fat] (-0.80902, -2.4899) -- (-0.49999, -3.4410) -- (-1.3090, -4.0287) -- (-1.6180, -3.0777) -- cycle;
\fill[fat] (2.1180, -1.5388) -- (3.1180, -1.5388) -- (3.4270, -2.4899) -- (2.4271, -2.4899) -- cycle;
\fill[fat] (2.1180, 1.5388) -- (2.4270, 2.4899) -- (3.4271, 2.4899) -- (3.1180, 1.5388) -- cycle;
\fill[fat] (-0.80902, 2.4899) -- (-0.50000, 3.4410) -- (-1.3090, 4.0287) -- (-1.6180, 3.0777) -- cycle;
\fill[fat] (-2.6180, -1.9021) -- (-3.4270, -2.4899) -- (-3.7361, -1.5388) -- (-2.9270, -0.95106) -- cycle;
\fill[fat] (-2.6180, -1.9021) -- (-3.4270, -2.4899) -- (-2.6180, -3.0777) -- (-1.8090, -2.4899) -- cycle;
\fill[fat] (1.0000, -3.0777) -- (1.3090, -4.0287) -- (0.30902, -4.0287) -- (1.9073e-6, -3.0777) -- cycle;
\fill[fat] (1.0000, -3.0777) -- (1.3090, -4.0287) -- (2.1180, -3.4410) -- (1.8090, -2.4899) -- cycle;
\fill[fat] (3.2361, 1.9073e-6) -- (4.2361, 0.00000) -- (3.9271, -0.95105) -- (2.9270, -0.95105) -- cycle;
\fill[fat] (3.2361, 1.9073e-6) -- (2.9270, 0.95106) -- (3.9270, 0.95106) -- (4.2361, 0.00000) -- cycle;
\fill[fat] (1.0000, 3.0777) -- (1.8090, 2.4899) -- (2.1180, 3.4410) -- (1.3090, 4.0287) -- cycle;
\fill[fat] (1.0000, 3.0777) -- (-1.9073e-6, 3.0777) -- (0.30902, 4.0287) -- (1.3090, 4.0287) -- cycle;
\fill[fat] (-2.6180, 1.9021) -- (-1.8090, 2.4899) -- (-2.6180, 3.0777) -- (-3.4271, 2.4899) -- cycle;
\fill[fat] (-2.6180, 1.9021) -- (-2.9271, 0.95105) -- (-3.7361, 1.5388) -- (-3.4271, 2.4899) -- cycle;
\fill[fat] (-3.4271, -0.58779) -- (-4.2361, -9.5367e-7) -- (-4.5451, -0.95106) -- (-3.7361, -1.5388) -- cycle;
\fill[fat] (-3.7361, -1.5388) -- (-3.4270, -2.4899) -- (-4.4271, -2.4899) -- (-4.7361, -1.5388) -- cycle;
\fill[fat] (-0.49999, -3.4410) -- (0.30902, -4.0287) -- (-0.50000, -4.6165) -- (-1.3090, -4.0287) -- cycle;
\fill[fat] (0.30902, -4.0287) -- (4.7684e-6, -4.9798) -- (1.0000, -4.9798) -- (1.3090, -4.0287) -- cycle;
\fill[fat] (3.1180, -1.5388) -- (3.4270, -2.4899) -- (4.2361, -1.9021) -- (3.9271, -0.95105) -- cycle;
\fill[fat] (3.9271, -0.95105) -- (4.7361, -1.5388) -- (5.0451, -0.58778) -- (4.2361, 0.00000) -- cycle;
\fill[fat] (2.4270, 2.4899) -- (3.4271, 2.4899) -- (3.1180, 3.4410) -- (2.1180, 3.4410) -- cycle;
\fill[fat] (2.1180, 3.4410) -- (2.9270, 4.0287) -- (2.1180, 4.6165) -- (1.3090, 4.0287) -- cycle;
\fill[fat] (-1.6180, 3.0777) -- (-1.3090, 4.0287) -- (-2.3090, 4.0287) -- (-2.6180, 3.0777) -- cycle;
\fill[fat] (-2.6180, 3.0777) -- (-2.9271, 4.0287) -- (-3.7361, 3.4409) -- (-3.4271, 2.4899) -- cycle;
\fill[fat] (-3.4270, -2.4899) -- (-3.7361, -3.4410) -- (-2.9270, -4.0287) -- (-2.6180, -3.0777) -- cycle;
\fill[fat] (1.3090, -4.0287) -- (2.1180, -4.6165) -- (2.9271, -4.0287) -- (2.1180, -3.4410) -- cycle;
\fill[fat] (4.2361, 0.00000) -- (5.0451, 0.58779) -- (4.7361, 1.5388) -- (3.9270, 0.95106) -- cycle;
\fill[fat] (1.3090, 4.0287) -- (0.30902, 4.0287) -- (-9.5367e-7, 4.9798) -- (1.0000, 4.9798) -- cycle;
\fill[fat] (-3.4271, 2.4899) -- (-4.4271, 2.4899) -- (-4.7361, 1.5388) -- (-3.7361, 1.5388) -- cycle;
\fill[fat] (-2.6180, -3.0777) -- (-1.6180, -3.0777) -- (-1.3090, -4.0287) -- (-2.3090, -4.0287) -- cycle;
\fill[fat] (2.1180, -3.4410) -- (2.4271, -2.4899) -- (3.4270, -2.4899) -- (3.1180, -3.4410) -- cycle;
\fill[fat] (3.9270, 0.95106) -- (3.1180, 1.5388) -- (3.4271, 2.4899) -- (4.2361, 1.9021) -- cycle;
\fill[fat] (0.30902, 4.0287) -- (-0.50000, 3.4410) -- (-1.3090, 4.0287) -- (-0.50000, 4.6165) -- cycle;
\fill[fat] (-3.7361, 1.5388) -- (-4.5451, 0.95105) -- (-4.2361, -9.5367e-7) -- (-3.4270, 0.58778) -- cycle;
\fill[fat] (-1.3090, -4.0287) -- (-2.3090, -4.0287) -- (-2.6180, -4.9798) -- (-1.6180, -4.9798) -- cycle;
\fill[fat] (-1.3090, -4.0287) -- (-0.50000, -4.6165) -- (-0.80901, -5.5676) -- (-1.6180, -4.9798) -- cycle;
\fill[fat] (3.4270, -2.4899) -- (4.2361, -3.0777) -- (3.9271, -4.0287) -- (3.1180, -3.4410) -- cycle;
\fill[fat] (3.4270, -2.4899) -- (4.2361, -3.0777) -- (5.0451, -2.4899) -- (4.2361, -1.9021) -- cycle;
\fill[fat] (3.4271, 2.4899) -- (4.2361, 3.0777) -- (3.9270, 4.0287) -- (3.1180, 3.4410) -- cycle;
\fill[fat] (3.4271, 2.4899) -- (4.2361, 3.0777) -- (5.0451, 2.4899) -- (4.2361, 1.9021) -- cycle;
\fill[fat] (-1.3090, 4.0287) -- (-2.3090, 4.0287) -- (-2.6180, 4.9798) -- (-1.6180, 4.9798) -- cycle;
\fill[fat] (-1.3090, 4.0287) -- (-1.6180, 4.9798) -- (-0.80902, 5.5676) -- (-0.50000, 4.6165) -- cycle;
\fill[fat] (-4.2361, -9.5367e-7) -- (-4.5451, 0.95105) -- (-5.5451, 0.95105) -- (-5.2361, -1.9073e-6) -- cycle;
\fill[fat] (-4.2361, -9.5367e-7) -- (-4.5451, -0.95106) -- (-5.5451, -0.95106) -- (-5.2361, -1.9073e-6) -- cycle;
	\node[joint] (0) at (0.00000, 0.00000) {};
	\node[joint] (1) at (-1.6180, -4.7684e-7) {};
	\node[joint] (2) at (-0.80902, -0.58779) {};
	\node[joint] (3) at (-0.50000, -1.5388) {};
	\node[joint] (4) at (0.30902, -0.95106) {};
	\node[joint] (5) at (1.3090, -0.95106) {};
	\node[joint] (6) at (1.0000, 0.00000) {};
	\node[joint] (7) at (1.3090, 0.95106) {};
	\node[joint] (8) at (0.30902, 0.95106) {};
	\node[joint] (9) at (-0.50000, 1.5388) {};
	\node[joint] (10) at (-0.80902, 0.58779) {};
	\node[joint] (11) at (-1.3090, -0.95106) {};
	\node[joint] (12) at (0.50000, -1.5388) {};
	\node[joint] (13) at (1.6180, 1.9073e-6) {};
	\node[joint] (14) at (0.50000, 1.5388) {};
	\node[joint] (15) at (-1.3090, 0.95106) {};
	\node[joint] (16) at (-2.6180, -9.5367e-7) {};
	\node[joint] (17) at (-2.3090, 0.95105) {};
	\node[joint] (18) at (-0.80902, -2.4899) {};
	\node[joint] (19) at (-1.6180, -1.9021) {};
	\node[joint] (20) at (2.1180, -1.5388) {};
	\node[joint] (21) at (1.3090, -2.1266) {};
	\node[joint] (22) at (2.1180, 1.5388) {};
	\node[joint] (23) at (2.4270, 0.58779) {};
	\node[joint] (24) at (-0.80902, 2.4899) {};
	\node[joint] (25) at (0.19098, 2.4899) {};
	\node[joint] (26) at (-2.3090, -0.95106) {};
	\node[joint] (27) at (0.19098, -2.4899) {};
	\node[joint] (28) at (2.4271, -0.58778) {};
	\node[joint] (29) at (1.3090, 2.1266) {};
	\node[joint] (30) at (-1.6180, 1.9021) {};
	\node[joint] (31) at (-2.6180, -1.9021) {};
	\node[joint] (32) at (1.0000, -3.0777) {};
	\node[joint] (33) at (3.2361, 1.9073e-6) {};
	\node[joint] (34) at (1.0000, 3.0777) {};
	\node[joint] (35) at (-2.6180, 1.9021) {};
	\node[joint] (36) at (-1.8090, -2.4899) {};
	\node[joint] (37) at (1.9073e-6, -3.0777) {};
	\node[joint] (38) at (1.8090, -2.4899) {};
	\node[joint] (39) at (2.9270, -0.95105) {};
	\node[joint] (40) at (2.9270, 0.95106) {};
	\node[joint] (41) at (1.8090, 2.4899) {};
	\node[joint] (42) at (-1.9073e-6, 3.0777) {};
	\node[joint] (43) at (-1.8090, 2.4899) {};
	\node[joint] (44) at (-2.9271, 0.95105) {};
	\node[joint] (45) at (-2.9270, -0.95106) {};
	\node[joint] (46) at (-3.4271, -0.58779) {};
	\node[joint] (47) at (-3.7361, -1.5388) {};
	\node[joint] (48) at (-0.49999, -3.4410) {};
	\node[joint] (49) at (0.30902, -4.0287) {};
	\node[joint] (50) at (3.1180, -1.5388) {};
	\node[joint] (51) at (3.9271, -0.95105) {};
	\node[joint] (52) at (2.4270, 2.4899) {};
	\node[joint] (53) at (2.1180, 3.4410) {};
	\node[joint] (54) at (-1.6180, 3.0777) {};
	\node[joint] (55) at (-2.6180, 3.0777) {};
	\node[joint] (56) at (-3.4270, -2.4899) {};
	\node[joint] (57) at (1.3090, -4.0287) {};
	\node[joint] (58) at (4.2361, 0.00000) {};
	\node[joint] (59) at (1.3090, 4.0287) {};
	\node[joint] (60) at (-3.4271, 2.4899) {};
	\node[joint] (61) at (-2.6180, -3.0777) {};
	\node[joint] (62) at (2.1180, -3.4410) {};
	\node[joint] (63) at (3.9270, 0.95106) {};
	\node[joint] (64) at (0.30902, 4.0287) {};
	\node[joint] (65) at (-3.7361, 1.5388) {};
	\node[joint] (66) at (-1.6180, -3.0777) {};
	\node[joint] (67) at (2.4271, -2.4899) {};
	\node[joint] (68) at (3.1180, 1.5388) {};
	\node[joint] (69) at (-0.50000, 3.4410) {};
	\node[joint] (70) at (-3.4270, 0.58778) {};
	\node[joint] (71) at (-1.3090, -4.0287) {};
	\node[joint] (72) at (3.4270, -2.4899) {};
	\node[joint] (73) at (3.4271, 2.4899) {};
	\node[joint] (74) at (-1.3090, 4.0287) {};
	\node[joint] (75) at (-4.2361, -9.5367e-7) {};
	\node[joint] (76) at (-2.3090, -4.0287) {};
	\node[joint] (77) at (3.1180, -3.4410) {};
	\node[joint] (78) at (4.2361, 1.9021) {};
	\node[joint] (79) at (-0.50000, 4.6165) {};
	\node[joint] (80) at (-4.5451, 0.95105) {};
	\node[joint] (81) at (-0.50000, -4.6165) {};
	\node[joint] (82) at (4.2361, -1.9021) {};
	\node[joint] (83) at (3.1180, 3.4410) {};
	\node[joint] (84) at (-2.3090, 4.0287) {};
	\node[joint] (85) at (-4.5451, -0.95106) {};
	\node[joint] (86) at (-1.6180, -4.9798) {};
	\node[joint] (87) at (-0.80901, -5.5676) {};
	\node[joint] (88) at (4.2361, -3.0777) {};
	\node[joint] (89) at (5.0451, -2.4899) {};
	\node[joint] (90) at (4.2361, 3.0777) {};
	\node[joint] (91) at (3.9270, 4.0287) {};
	\node[joint] (92) at (-1.6180, 4.9798) {};
	\node[joint] (93) at (-2.6180, 4.9798) {};
	\node[joint] (94) at (-5.2361, -1.9073e-6) {};
	\node[joint] (95) at (-5.5451, -0.95106) {};
	\node[joint] (96) at (4.7684e-6, -4.9798) {};
	\node[joint] (97) at (4.7361, -1.5388) {};
	\node[joint] (98) at (2.9270, 4.0287) {};
	\node[joint] (99) at (-2.9271, 4.0287) {};
	\node[joint] (100) at (-4.7361, -1.5388) {};
	\node[joint] (101) at (1.0000, -4.9798) {};
	\node[joint] (102) at (5.0451, -0.58778) {};
	\node[joint] (103) at (2.1180, 4.6165) {};
	\node[joint] (104) at (-3.7361, 3.4409) {};
	\node[joint] (105) at (-4.4271, -2.4899) {};
	\node[joint] (106) at (-2.6180, -4.9798) {};
	\node[joint] (107) at (3.9271, -4.0287) {};
	\node[joint] (108) at (5.0451, 2.4899) {};
	\node[joint] (109) at (-0.80902, 5.5676) {};
	\node[joint] (110) at (-5.5451, 0.95105) {};
	\node[joint] (111) at (-2.9270, -4.0287) {};
	\node[joint] (112) at (2.9271, -4.0287) {};
	\node[joint] (113) at (4.7361, 1.5388) {};
	\node[joint] (114) at (-9.5367e-7, 4.9798) {};
	\node[joint] (115) at (-4.7361, 1.5388) {};
	\node[joint] (116) at (-3.7361, -3.4410) {};
	\node[joint] (117) at (2.1180, -4.6165) {};
	\node[joint] (118) at (5.0451, 0.58779) {};
	\node[joint] (119) at (1.0000, 4.9798) {};
	\node[joint] (120) at (-4.4271, 2.4899) {};
	\node[joint] (121) at (-1.8090, -5.5676) {};
	\node[joint] (122) at (4.7361, -3.4410) {};
	\node[joint] (123) at (4.7361, 3.4410) {};
	\node[joint] (124) at (-1.8090, 5.5676) {};
	\node[joint] (125) at (-5.8541, -9.5367e-7) {};
	\draw[edge](0)edge(2) (0)edge(4) (0)edge(6) (0)edge(8) (0)edge(10) (1)edge(2) (1)edge(10) (1)edge(11) (1)edge(15) (1)edge(16) (2)edge(3) 
	(3)edge(4) (3)edge(11) (3)edge(12) (3)edge(18) (4)edge(5) (5)edge(6) (5)edge(12) (5)edge(13) (5)edge(20) (6)edge(7) (7)edge(8) (7)edge(13)
	 (7)edge(14) (7)edge(22) (8)edge(9) (9)edge(10) (9)edge(14) (9)edge(15) (9)edge(24) (11)edge(19) (11)edge(26) (12)edge(21) (12)edge(27) 
	 (13)edge(23) (13)edge(28) (14)edge(25) (14)edge(29) (15)edge(17) (15)edge(30) (16)edge(17) (16)edge(26) (16)edge(44) (16)edge(45) 
	 (16)edge(46) (16)edge(70) (17)edge(35) (18)edge(19) (18)edge(27) (18)edge(36) (18)edge(37) (18)edge(48) (18)edge(66) (19)edge(31) 
	 (20)edge(21) (20)edge(28) (20)edge(38) (20)edge(39) (20)edge(50) (20)edge(67) (21)edge(32) (22)edge(23) (22)edge(29) (22)edge(40) 
	 (22)edge(41) (22)edge(52) (22)edge(68) (23)edge(33) (24)edge(25) (24)edge(30) (24)edge(42) (24)edge(43) (24)edge(54) (24)edge(69) 
	 (25)edge(34) (26)edge(31) (27)edge(32) (28)edge(33) (29)edge(34) (30)edge(35) (31)edge(36) (31)edge(45) (31)edge(56) (32)edge(37) 
	 (32)edge(38) (32)edge(57) (33)edge(39) (33)edge(40) (33)edge(58) (34)edge(41) (34)edge(42) (34)edge(59) (35)edge(43) (35)edge(44) 
	 (35)edge(60) (36)edge(61) (37)edge(49) (38)edge(62) (39)edge(51) (40)edge(63) (41)edge(53) (42)edge(64) (43)edge(55) (44)edge(65) 
	 (45)edge(47) (46)edge(47) (46)edge(75) (47)edge(56) (47)edge(85) (47)edge(100) (48)edge(49) (48)edge(71) (49)edge(57) (49)edge(81) 
	 (49)edge(96) (50)edge(51) (50)edge(72) (51)edge(58) (51)edge(82) (51)edge(97) (52)edge(53) (52)edge(73) (53)edge(59) (53)edge(83) 
	 (53)edge(98) (54)edge(55) (54)edge(74) (55)edge(60) (55)edge(84) (55)edge(99) (56)edge(61) (56)edge(105) (56)edge(116) (57)edge(62)
	  (57)edge(101) (57)edge(117) (58)edge(63) (58)edge(102) (58)edge(118) (59)edge(64) (59)edge(103) (59)edge(119) (60)edge(65) 
	  (60)edge(104) (60)edge(120) (61)edge(66) (61)edge(76) (61)edge(111) (62)edge(67) (62)edge(77) (62)edge(112) (63)edge(68) 
	  (63)edge(78) (63)edge(113) (64)edge(69) (64)edge(79) (64)edge(114) (65)edge(70) (65)edge(80) (65)edge(115) (66)edge(71) 
	  (67)edge(72) (68)edge(73) (69)edge(74) (70)edge(75) (71)edge(76) (71)edge(81) (71)edge(86) (72)edge(77) (72)edge(82) (72)edge(88) 
	  (73)edge(78) (73)edge(83) (73)edge(90) (74)edge(79) (74)edge(84) (74)edge(92) (75)edge(80) (75)edge(85) (75)edge(94) (76)edge(106)
	   (77)edge(107) (78)edge(108) (79)edge(109) (80)edge(110) (81)edge(87) (82)edge(89) (83)edge(91) (84)edge(93) (85)edge(95) 
	   (86)edge(87) (86)edge(106) (87)edge(96) (87)edge(121) (88)edge(89) (88)edge(107) (89)edge(97) (89)edge(122) (90)edge(91) 
	   (90)edge(108) (91)edge(98) (91)edge(123) (92)edge(93) (92)edge(109) (93)edge(99) (93)edge(124) (94)edge(95) (94)edge(110) 
	   (95)edge(100) (95)edge(125) (96)edge(101) (97)edge(102) (98)edge(103) (99)edge(104) (100)edge(105) (106)edge(111) (106)edge(121)
	    (107)edge(112) (107)edge(122) (108)edge(113) (108)edge(123) (109)edge(114) (109)edge(124) (110)edge(115) (110)edge(125) 
	    (111)edge(116) (112)edge(117) (113)edge(118) (114)edge(119) (115)edge(120);
	\end{tikzpicture}
	\caption{A piece of a Penrose tiling.}
	\label{fig:Penrosetiling}
\end{figure}

\begin{figure}[ht]
\begin{center}
	\begin{tikzpicture}[scale=2]
		\node[fvertex] (1T) at (0,0.588) {};
		\node[fvertex] (2T) at (0.809,0) {};
		\node[fvertex] (3T) at (0,-0.588) {};
		\node[fvertex] (4T) at (-0.809,0) {};
		\node at (0,0) {$F$};
		
		\draw[edge,dashed,->] (1T) -- (2T);
		\draw[edge,dashed,<-] (2T) -- (3T);
		\draw[edge,<-] (3T) -- (4T);
		\draw[edge,->] (4T) -- (1T);
		\node[fvertex] (1t) at (3,0.309) {};
		\node[fvertex] (2t) at (3.951,0) {};
		\node[fvertex] (3t) at (3,-0.309) {};
		\node[fvertex] (4t) at (2.049,0) {};
		\node at (3,0) {$T$};
		
		\draw[edge,dashed,<-] (1t) -- (2t);
		\draw[edge,->] (2t) -- (3t);
		\draw[edge,<-] (3t) -- (4t);
		\draw[edge,dashed,->] (4t) -- (1t);
	\end{tikzpicture}
\end{center}
	\caption{The two types of rhombi used in a Penrose tiling with their matching rules;
	the ``fat'' rhombus $F$ (left) and the ``thin'' rhombus $T$ (right).}\label{rhombus}
\end{figure}
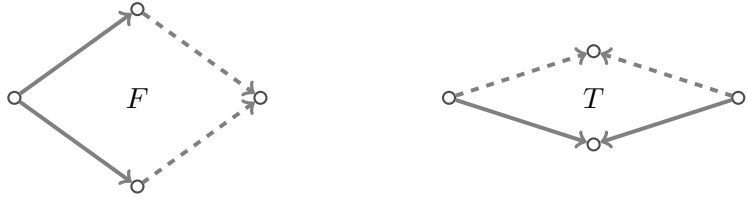

We shall apply the results of \Cref{sec:braced} to the family of infinite P-frameworks generated by the Penrose (rhombus) tilings; see \Cref{fig:Penrosetiling}.
Any Penrose tiling is formed from two types of rhombi with sides marked by arrows (see Figure \ref{rhombus}):
a ``fat'' rhombus tile $F$ with angles $72^\circ$ and $108^\circ$,
and a ``thin'' rhombus tile $T$ with angles $36^\circ$ and $144^\circ$.
To construct a Penrose tiling, one simply needs to match copies of the tiles together so that the arrows they meet at are oriented in the same direction and have the same type (i.e.~dashed or not dashed).
The copies of the two types of tiles can be rotated, reflected and translated however you like.
For more details on the described construction, we refer the reader to~\cite{DEBRUIJN198139}.
Interestingly, while there exist uncountably many distinct tilings (modulo rotations, reflections and translations) that can be constructed by this method,
there exist exactly two which have 5-fold rotational symmetry.
Both of these properties can be proven by de Brujin’s method of implementing \emph{pentagrids};
see \cite{DEBRUIJN198139,DEBRUIJN198153} for more information on the topic.

By considering the shared sides of tiles in a Penrose tiling as edges and the shared corners as vertices,
it is easy to see how we can transform any Penrose tiling into an infinite P-framework.
Any such P-framework that we obtain from this method shall be called a \emph{Penrose framework}.

Let $(G, \rho)$ be a P-framework with ribbon $r$.
As every edge in a ribbon is parallel,
we can choose a unit vector $u_r$ so that $u_r$ is orthogonal to $\rho(u)-\rho(v)$ for all $uv \in r$;
we call any such unit vector a \emph{direction} of~$r$.
It is immediate that every ribbon has exactly two directions: $u_r$ and $-u_r$.
We say two ribbons $r,r'$ with directions $u_r,u_{r'}$ respectively are \emph{parallel} if $u_r$ and $u_{r'}$ are parallel.
The following result can be proven by de Brujin's method of implementing pentagrids \cite{DEBRUIJN198139,DEBRUIJN198153}.

\begin{lemma}\label{l:penroseribbons}
	Let $(G,\rho)$ be a Penrose framework with ribbon graph $\Gamma$.
	Then the following properties hold.
	\begin{enumerate}[(i)]
		\item\label{l:penroseribbons1} Two ribbons are not adjacent in $\Gamma$ if and only if they are parallel.
		\item\label{l:penroseribbons2}	There exist 5 cyclically ordered unit vectors $u_1,\ldots,u_5$ where the angle between $u_{i \bmod 5}$ and $u_{i+1 \bmod 5}$ is $2\pi/5$,
		so that any ribbon $r$ has direction $u_i$ for some $i$.
		\item\label{l:penroseribbons2.5} For each unit vector $u_i$, there exists a countably infinite amount of ribbons with direction $u_i$.
		\item\label{l:penroseribbons3} If $r$ is a ribbon with direction $u_{i \bmod 5}$ and $r'$ is a ribbon with direction $u_{i+1 \bmod 5}$,
		then there exists a unique 4-cycle in $r \cup r'$,
		and its image under $\rho$ is a copy of the tile $F$.
		\item\label{l:penroseribbons4} If $r$ is a ribbon with direction $u_{i \bmod 5}$ and $r'$ is a ribbon with direction $u_{i+2 \bmod 5}$,
		then there exists a unique 4-cycle in $r \cup r'$,
		and its image under $\rho$ is a copy of the tile $T$.
	\end{enumerate}
\end{lemma}

Similar to the infinite braced grid, we can achieve rigidity by simply bracing every parallelogram in two ribbons that are adjacent in the ribbon graph.

\begin{corollary}\label{cor:penroseRibbon}
	Let $(G,\rho)$ be a braced Penrose framework formed by bracing every rhombus in any two ribbons.
	Then $(G,\rho)$ is flexible if and only if the two ribbons are parallel.
\end{corollary}

\begin{proof}
	Let $r_1,r_2$ be the braced ribbons and let $\Gamma_b$ be the bracing graph of $(G, \rho)$.
	By \Cref{l:penroseribbons}, each ribbon $r_i$ is adjacent in $\Gamma_b$ to every ribbon it is not parallel to.
	If $r_1,r_2$ are not parallel, then $r_1 r_2$ is an edge of $\Gamma_b$ and every other vertex in $\Gamma_b$ is adjacent to $r_1$ or~$r_2$,
	hence $(G,\rho)$ is rigid by \Cref{thm:bracedPframeworks}.
	If $r_1,r_2$ are parallel, then every ribbon parallel to both $r_1,r_2$ will be an isolated vertex in $\Gamma_b$,
	hence $(G,\rho)$ is flexible by \Cref{thm:bracedPframeworks}.
\end{proof}

The following result was proven by Power for infinitesimal rigidity \cite{power2021parallelogram}.
It follows then that we can achieve the following result by combining Power's result with \Cref{c:parallelogram},
however we have included a proof for the sake of completeness.

\begin{corollary}\label{cor:penrose}
	Let $(G,\rho)$ be a braced Penrose framework formed by bracing every rhombus that corresponds to the tile $F$,
	or by bracing every rhombus that corresponds to the tile $T$.
	Then $(G,\rho)$ is rigid.
\end{corollary}

\begin{proof}
	Define the graph $H = (V_H,E_H)$, where $V_H$ is the disjoint union of countably infinite sets $V_1,\ldots,V_5$, and $uv \in E_H$ for $u \in V_i$, $v \in V_j$ if and only if $(i - j) \equiv \pm 1 \bmod 5$.
	By \Cref{l:penroseribbons},
	the bracing graph of $(G,\rho)$ is isomorphic to $H$, irrespective of whether 
	we are bracing every rhombus corresponding to the tile~$F$ or to the tile $T$.
	As the graph~$H$ is connected it follows that $(G,\rho)$ is rigid by \Cref{thm:bracedPframeworks}.
\end{proof}

In fact, as we shall see with the following two results,
we do not need to brace all the rhombi corresponding to a given type of tile to obtain a rigid braced Penrose framework.

\begin{corollary}\label{cor:penrose2}
	Let $P$ be a Penrose tiling.
	Choose a tile type $S \in \{F,T\}$ and a rotational orientation $\sigma$ so that $P$ contains a translate of the tile $\sigma S$.
	Let $(G,\rho)$ be a braced Penrose framework formed by bracing every rhombus that corresponds to the tile $S$ except those that are translates of $\sigma S$.
	Then $(G,\rho)$ is rigid.
\end{corollary}

\begin{proof}
	Define the graph $H = (V_H,E_H)$, where $V_H$ is the disjoint union of countably infinite sets $V_1,\ldots,V_5$, and $uv \in E_H$ for $u \in V_i$, $v \in V_j$ if and only if $|i-j| = 1$.
	By \Cref{l:penroseribbons},
	the bracing graph of $(G,\rho)$ is isomorphic to $H$.
	As the graph $H$ is connected it follows that $(G,\rho)$ is rigid by \Cref{thm:bracedPframeworks}.
\end{proof}

\begin{corollary}\label{cor:penroseProb}
	Fix a tile $S \in \{F,T\}$. Let $(G_p,\rho)$ be a braced Penrose framework formed by randomly bracing every rhombus
	that corresponds to the tile $S$ with probability $p \in (0,1]$.
	Then $(G_p,\rho)$ is almost surely rigid.
\end{corollary}

\begin{proof}
	Define the graph $H = (V_H,E_H)$, where $V_H$ is the disjoint union of countably infinite sets $V_1,\ldots,V_5$, and $uv \in E_H$ for $u \in V_i$, $v \in V_j$ if and only if $(i - j) \equiv \pm 1 \bmod 5$.
	By \Cref{l:penroseribbons},
	the bracing graph of $(G_1,\rho)$ is isomorphic to $H$.
	Now define $H_p$ be a spanning subgraph of $H$ that corresponds to the bracing graph of $(G_p , \rho)$ under the isomorphism between the ribbon graph and $H$.
	It follows from \Cref{thm:bracedPframeworks} that we need only show that $H_p$ is connected.
	Note that the probability that an edge of $H$ is also an edge of $H_p$ is $p$.
	
	Choose two vertices $u,v \in V_H$.
	We may assume by relabelling if necessary that $u \in V_1$ and $v \in V_t$ for some $t \in \{1,2,3\}$.
	First suppose that $t \in \{1,3\}$.
	The probability of a vertex $w \in V_2$ being adjacent to both $u$ and $v$ is $p^2$,
	hence the probability that there exists a path from $u$ to $v$ in $H_p$ is at least $\lim_{n \rightarrow \infty} 1 - (1-p^2)^n = 1$.
	Now suppose that~$t= 2$.
	Given vertices $w_1 \in V_1$ and $w_2 \in V_2$,
	the probability that $(u, w_2 , w_1, v)$ is a path in $H_p$ is $p^3$,
	hence the probability that there exists a path from $u$ to $v$ in $H_p$ is at least $\lim_{n \rightarrow \infty} 1 - (1-p^3)^n = 1$.
\end{proof}

\section{Rotationally symmetric infinite frameworks}\label{sec:rot}
In this section, we extend the main result of~\cite{RotSymmetry} to infinite graphs.
We recall the definitions of $k$-fold rotational symmetry, where $k \geq 2$.

\begin{definition}
Let $G$ be a graph and let the group $\Ck := \left\langle \omega : \omega^k =1 \right\rangle$
act on $G$, namely, there exists an injective group homomorphism $\varphi: \Ck \rightarrow \Aut(G)$.
We define $\gamma v := \varphi(\gamma)(v)$ for $\gamma\in\Ck$ 
and $\gamma e := \gamma u \gamma v$ for any edge $e=uv\in E_G$.
A vertex $v \in V_G$ is an \emph{invariant vertex} if $\gamma v =v$ for all $\gamma \in \Ck$,
and \emph{partially invariant} if $\gamma v =v$ for some $\gamma \in \Ck$, $\gamma\neq 1$.
The graph $G$ is called \emph{$\Ck$-symmetric} if every partially invariant vertex is invariant,
and the set of invariant vertices of $G$ forms an independent set.
\end{definition}

\begin{definition}
A realization $\rho$ of a $\Ck$-symmetric graph $G$ in $\RR^2$ is called \emph{$\Ck$-symmetric} 
if $\rho(\gamma v ) = \Rot_\gamma\rho(v)$ for each $v \in V_G$ and $\gamma = \omega^j\in\Ck$,
where $\Rot_\omega$ is the $2\pi/k$ clockwise rotation matrix and $\Rot_\gamma = \Rot_{\omega}^j$.
The framework $(G,\rho)$ is \emph{$\Ck$-symmetric} if $\rho$ is $\Ck$-symmetric. 
\end{definition}

Notice that the requirements on invariant vertices to form an independent set is
justified by the fact that they must be mapped to the origin in a $\Ck$-symmetric framework
and we disallow edges to have zero length.

\begin{definition}
If there is a non-trivial flex $\rho_t$ of a $\Ck$-symmetric framework $(G,\rho)$
such that each $\rho_t$ is $\Ck$-symmetric,
then $(G,\rho)$ is \emph{$\Ck$-symmetric flexible} (or \emph{$k$-fold rotation symmetric flexible}),
and \emph{$\Ck$-symmetric rigid} otherwise.
\end{definition}

\begin{definition}\label{defn:NACnfold}
	Let $\delta$ be a NAC-coloring of a $\Ck$-symmetric graph $G$.
	A \emph{\red{} component}~$H$, i.e., a connected component of $G_\red^\delta$,
	is \emph{partially invariant} if there exists $\gamma \in \Ck \setminus \{1\}$ such that $\gamma H = H$,
	and \emph{invariant} if $\gamma H = H$ for all $\gamma \in \Ck$.
	We define \blue{} partially invariant components and \blue{} invariant components analogously.
	The NAC-coloring $\delta$ is called \emph{$\Ck$-symmetric} if
	$\delta(\gamma e) = \delta(e)$ for all $e \in E_G$ and $\gamma \in \Ck$ and
	no two distinct \blue{}, resp.\ \red{}, partially invariant components are connected by an edge.
\end{definition}
Similarly to standard NAC-colorings, $\Ck$-symmetric NAC-colorings can be used to determine
whether a $\Ck$-symmetric graph has a $\Ck$-symmetric flexible realization.
 
\begin{theorem}[\cite{RotSymmetry}]\label{thm:finitesymm}
	A finite $\Ck$-symmetric connected graph has a $\Ck$-symmetric
	NAC-coloring if and only if it has a $\Ck$-symmetric flexible realization in $\RR^2$. 
\end{theorem}
 
We shall now extend this result to infinite $\Ck$-symmetric graphs.
We begin with the following construction that was originally outlined in \cite{RotSymmetry} for finite $\Ck$-symmetric graphs.

\begin{proposition}
\label{prop:CnConstruction}
	Let $G$ be a countably infinite $\Ck$-symmetric graph.
	If $G$ has a $\Ck$-symmetric NAC-coloring, then there is a $\Ck$-symmetrically flexible realization of $G$.
\end{proposition}

\begin{proof}
	Let $R^0_1, \dots, R^{k-1}_1, R^0_2, \dots, R^{k-1}_2, \dots$ (resp.\ $B^0_1, \dots, B^{k-1}_1,\dots$)
	be the \red{} (resp.\ \blue{}) components of $G_\red^\delta$ (resp.\ $G_\blue^\delta$)
	that are not partially invariant.
	We can assume that $R^i_j = \omega^i R^0_j$ and $B^i_j = \omega^i B^0_j$ for $0\leq i < k$ and $1\leq j$.
	
	Let $a_1, a_2, \dots$ and $b_1,b_2,\dots$ be points in $\RR^2\setminus\{(0,0)\}$ 
	such that $a_j \neq \Rot_\omega^i a_{j'}$ and
	$b_j \neq \Rot_\omega^i b_{j'}$
	for $j\neq j'$ and $0\leq i<k$ arbitrary.
	We define functions $\overline{a},\overline{b} \colon V_G\rightarrow \RR^2$ by
	\begin{equation*}
		\overline{a}(v) = \begin{cases}
			\Rot_\omega^i a_j &\text{if } v \in R^i_j \\
			(0,0) 		 &\text{otherwise,}
		\end{cases}
	\quad\text{ and }\quad
	 \overline{b}(v) = \begin{cases}
			\Rot_\omega^i b_j &\text{if } v \in B^i_j \\
			(0,0) 		 &\text{otherwise.}
		\end{cases}
	\end{equation*}
	Notice that a vertex is mapped to $(0,0)$
	if and only if it is in a partially invariant component.
	We now construct realizations $\rho_t$ of $G$ by setting
	\begin{equation*}
		p_t(v) := \Rot(t) \overline{a}(v) + \overline{b}(v) \, ,
	\end{equation*}
	where $\Rot(t)$ is the clockwise rotation matrix by $t$ radians.
	By the same arguments as in~\cite[Lemma~2]{RotSymmetry},
	we can see that this gives a $\Ck$-symmetric flex. 
\end{proof}

We will now need to define a special type of subgraph tower.

\begin{definition}
	Let $G$ be a $\Ck$-symmetric graph with a subgraph tower $(G_n)_{n \in \mathbb{N}}$.
	We say the subgraph tower is \emph{$\Ck$-symmetric} if for every subgraph $G_n$,
	the restriction of every automorphism $\varphi(\gamma)$ in the image of $\varphi$ to $G_n$ is also an automorphism of $G_n$.
\end{definition}

\begin{lemma}
	\label{lem:CnNACfromSequence}
	Let $G$ be a countably infinite $\Ck$-symmetric graph with a $\Ck$-symmetric subgraph tower $(G_n)_{n \in \NN}$.
	Let $e_1,e_2$ be two edges of $G_1$.
	If there is a $\Ck$-symmetric NAC-coloring $\delta_n$ of $G_n$ such that $\delta_n(e_1)\neq\delta_n(e_2)$ for all $n\in \NN$,
	then $G$ has a $\Ck$-symmetric NAC-coloring.	
\end{lemma}

\begin{proof}
	Let  
	\begin{equation*}
		S_n = \{\delta\in\nac{G_n} : \text{ $\delta$ is $\Ck$-symmetric and } \delta(e_1)=\blue \neq \delta(e_2)\} 
	\end{equation*}
	for every $n \in \NN$.
	Each $S_n$ is non-empty since either $\delta_n$ is in $S_n$,
	or the NAC-coloring obtained by swapping the colors from $\delta_n$ is in $S_n$. 
	Let $\prec$ be the relation in $\bigcup_{n \in \NN} S_n$ given by restriction to edges of the smaller graph.
	Clearly, a restriction of a NAC-coloring that is invariant under the $\Ck$ group action to a subgraph is a NAC-coloring that is invariant under the $\Ck$ group action,
	provided that the obtained coloring is surjective,
	which is guaranteed by the condition $\delta(e_1)=\blue \neq \delta(e_2)$.
	Notice that a monochromatic component $H$ is partially invariant if and only if
	the orbit of each non-invariant vertex in $H$ contains at least two elements from $V_H$.
	Since the restriction removes all vertices in the same orbit,
	the restriction of a monochromatic component cannot be partially invariant unless the unrestricted component is partially invariant.
	This guarantees that no two partially invariant monochromatic components in the restriction are connected by an edge,
	namely, the NAC-coloring obtained by the restriction is $\Ck$-symmetric.
	Hence, the assumption of \Cref{lemma:koenig} holds.
	Therefore, we get a sequence of $\Ck$-symmetric NAC-colorings of $(G_n)_{n \in \NN}$ extending each other.
	
	This sequence of $\Ck$-symmetric NAC-colorings gives a NAC-coloring for~$G$
	by the same reasons as in \Cref{lem:infNACfromSequence}.
	It is immediate that the NAC-coloring will also be invariant under the $\Ck$ group action.
	Any such NAC-coloring will also be $\Ck$-symmetric,
	since for every two vertices $u,v$ contained in partially invariant monochromatic components,
	we can restrict to a $G_n$ containing both $u$ and $v$ to see that $u$ and $v$ are non-adjacent.
\end{proof}

\begin{theorem}
	A $\Ck$-symmetric connected countably infinite graph has a $\Ck$-symmetric
	NAC-coloring if and only if it has $\Ck$-symmetric flexible realization in $\RR^2$. 
\end{theorem}

\begin{proof}
	If $G$ has a $\Ck$-symmetric NAC-coloring then there is a $\Ck$-symmetric flexible realization of $G$ by \Cref{prop:CnConstruction}.
	Now suppose $G$ has a $\Ck$-symmetric realization with a non-trivial $\Ck$-symmetric flex $\alpha$.
	Since $\alpha$ is non-trivial, there are edges $e_1,e_2 \in E_G$ whose angle changes along the flex.
	We choose any $\Ck$-symmetric subgraph tower $(G_n)_{n \in \mathbb{N}}$ such that $e_1,e_2 \in E_{G_n}$.
	Clearly, the restriction of the flex $\alpha$ to $V_{G_n}$ is a $\Ck$-symmetric flex of $G_n$
	such that the angle between $e_1$ and $e_2$ changes.   
	Hence by \Cref{lem:active} with the restricting set given by $y_{\bar{u}} =0$ for some non-invariant vertex $\bar{u} \in V_G$ and
	\begin{equation} \label{eq:symmetric}
		\begin{aligned}
			x_{\omega v}=  \cos(2\pi/k) x_v + \sin(2\pi/k) y_v   \,, \quad y_{\omega v}=  -\sin (2\pi/k) x_v + \cos(2\pi/k) y_v\,
		\end{aligned}
	\end{equation}
	for every $v \in V_G$,
	there is a NAC-coloring $\delta_n$ of $G_n$ such that $\delta_n(e_1)\neq\delta_n(e_2)$ for all~$n\in \NN$.
	By the technique of \cite[Lemma~1]{RotSymmetry}, every $\delta_n$ is $\Ck$-symmetric.
	A NAC-coloring of $G$ is now provided by \Cref{lem:CnNACfromSequence}.
\end{proof}

\section{Rotationally symmetric P-frameworks}\label{sec:rot+braced}

We now wish to combine the results of \Cref{sec:braced,sec:rot}.
There are, however, some technical results we need to first cover.
We begin with the two following lemmas which describe how ribbon edge cuts behave.
The results can seen to be direct extensions of Lemma~3.2 and Lemma~3.6 of \cite{Grasegger2020} to the setting of infinite P-frameworks. 

\begin{lemma}\label{lem:RibbonEdgeCut}
	Let $(G,\rho)$ be a (countably infinite) P-framework.
	If $r$ is a ribbon of $G$,
	then $G \setminus r$ has exactly two connected components.
\end{lemma}

\begin{lemma}\label{lem:RibbonTransl}
	Let $(G,\rho)$ be a (countably infinite) P-framework with a ribbon $r$.
	Let $V_1 \cup V_2$ be the vertex set of $r$, 
	where all the vertices of $V_i$ lie in the same connected component of $G \setminus r$.
	Then for every pair of edges $u_1u_2,v_1v_2 \in r$ with $u_i,v_i \in V_i$,
	we have $\rho(u_1)-\rho(u_2) = \rho(v_1)-\rho(v_2)$.
\end{lemma}

Our next key result tells us that $\Ck$-symmetric NAC-colorings can be defined much more simply for the graphs of P-frameworks.

\begin{lemma}\label{lem:PartInvCompsPframework}
	Let $(G, \rho)$ be a $\Ck$-symmetric connected (countably infinite) P-framework,
	and let $\delta$ be a cartesian NAC-coloring of $G$.
	If $\delta(\gamma e) = \delta(e)$ for all $e \in E_G$ and $\gamma \in \Ck$,
	then $\delta$ is a $\Ck$-symmetric NAC-coloring of $G$.
\end{lemma}

\begin{proof}
	Suppose without loss of generality that there exists two red partially invariant components $A_1$ and $A_2$
	which are connected by a blue edge $u_1 u_2$.
	Define $r$ to be the ribbon containing $u_1 u_2$.
	Define $\gamma_i$ to be a group element that $A_i$ is invariant under.
	Since $u_i$ and $\gamma_i u_i$ lie in the same component,
	there exists a path $P_i \subset E_G$ in the induced subgraph on $A_i$ from $u_i$ to $\gamma_i u_i$.
	We note that for every integer $j$,
	we have $\gamma_2(\gamma_1 \gamma_2)^{j-1} A_1=\gamma_2^jA_1$ connected to $(\gamma_1 \gamma_2)^j A_2=\gamma_1^jA_2$
	by the edge $((\gamma_1 \gamma_2)^j u_1)((\gamma_1 \gamma_2)^j u_2)$,
	and we have $(\gamma_1 \gamma_2)^j A_2=\gamma_1^jA_2$ connected to $\gamma_2(\gamma_1 \gamma_2)^j A_1=\gamma_2^{j+1}A_1$
	by the edge $(\gamma_2(\gamma_1 \gamma_2)^ju_1)(\gamma_2(\gamma_1 \gamma_2)^j u_2)$.
	As the order of $\gamma_1 \gamma_2$ is finite,
	there exists a circuit $C$ given by the sequence
	\begin{multline*}
		u_1 u_2, \, P_2, \,(\gamma_2 u_2) (\gamma_2 u_1), \, \gamma_2 P_1, \, (\gamma_1\gamma_2 u_1)(\gamma_1\gamma_2 u_2), \, \gamma_1 \gamma_2 P_2 \ldots \\
		\ldots, \,(\gamma_2(\gamma_1 \gamma_2)^{m-1}u_2) (\gamma_2(\gamma_1 \gamma_2)^{m-1}u_1), \,  \gamma_2(\gamma_1 \gamma_2)^{m-1} P_1, \, u_1 u_2
	\end{multline*}
	for some positive integer $m$.
	By \Cref{lem:RibbonEdgeCut},
	$G \setminus r$ has exactly two connected components,
	one containing $u_1$ and the other containing $u_2$.
	Because $u_1u_2$ is contained in the circuit~$C$,
	the ribbon $r$ must contain an edge of $C \setminus \{u_1u_2\}$.
	As each edge in each path $\gamma_2(\gamma_1\gamma_2)^j P_1$ and $(\gamma_1\gamma_2)^\ell P_2$ is red and each edge of $r$ is blue (\Cref{lem:cartesianIffRibbonsMonochromatic}),
	the ribbon $r$ must contain an edge
	$(\gamma u_1)(\gamma u_2)$ for some $\gamma \in \left\langle \gamma_1, \gamma_2 \right\rangle \setminus \{1\}$.
	As distinct ribbons do not intersect each other, we have $\gamma r = r$.
	The only possible two symmetries that a ribbon of a P-framework can have are $\C_2$-symmetric rotational symmetry and reflectional symmetry.
	Hence $\Rot_\gamma$ is the $\pi$ rotation matrix,
	i.e.~$\gamma^2=1$ and $\gamma \neq 1$.
	
	If both $\gamma_1$ and $\gamma_2$ have odd orders then, since $\Ck$ is a cyclic group,
	we have $\left\langle \gamma_1, \gamma_2 \right\rangle$ has an odd order\footnote{As $\Ck$ is cyclic, 
	there exists $\ell$ such that $\left\langle \gamma_1, \gamma_2 \right\rangle = \left\langle \omega^{\ell} \right\rangle$.
	If we choose $\ell$ to be the smallest positive integer with this property, then $k/\ell$ is the order of $\left\langle \gamma_1, \gamma_2 \right\rangle$. 
	Given $k_i$ is the smallest positive integer so that $\left\langle \gamma_i \right\rangle = \left\langle \omega^{k_i} \right\rangle$ for each $i$,
	we note that $k/k_i$ is the order of $\left\langle \gamma_i \right\rangle$ and, by B\'ezout's identity, $\ell = \gcd(k_1,k_2)$. 
	If we denote $k = a 2^b$ for some odd $a$ then, as $k_1,k_2$ have odd order, 
	we must have $k_1 = a_1 2^b$ and $k_2 = a_2 2^b$ for odd $a_1,a_2$, since both $k/k_1$ and $k/k_2$ are odd.
	Hence $\ell = c 2^b$ for some odd $c$ and $k/\ell$ is odd.},
	contradicting that $\gamma \in \left\langle \gamma_1, \gamma_2 \right\rangle$.
	Without loss of generality,
	we may assume the order of $\gamma_1$ is even.
	Hence,  $\left\langle \gamma_1 \right\rangle$ contains an element which is its own inverse, namely, it corresponds to rotation by $\pi$.
	Hence $\gamma \in \left\langle \gamma_1 \right\rangle$ and $\gamma A_1 = A_1$;
	in particular, $u_1, \gamma u_1 \in A_1$.
	Since $\rho( \gamma u_1) - \rho(\gamma u_2) = \rho(u_2) - \rho(u_1)$,
	we have by \Cref{lem:RibbonTransl} that $u_1$ and $\gamma u_1$ lie in separate connected components of $G \setminus r$.
	However, this contradicts that every pair of vertices in $A_1$ are connected by a path in $G \setminus r$;
	this follows from noting that no edge of $r$ can lie in $A_1$,
	since every edge in $r$ is blue by \Cref{lem:cartesianIffRibbonsMonochromatic}.
\end{proof}

The most technical part of this section is the construction of a $\Ck$-symmetric flex.
\begin{lemma}
	\label{lem:CnCartNACconstruction}
	If a braced $\Ck$-symmetric (countably infinite) P-framework has a cartesian $\Ck$-symmetric NAC-coloring,
	then it is $\Ck$-symmetrically flexible.
\end{lemma}
\begin{proof}
	We extend the construction introduced in \cite{Grasegger2020} to countably infinite frameworks
	and show that it yields a $\Ck$-symmetric flex for both finite and infinite case.
	Let $(G,\rho)$ be a braced $\Ck$-symmetric (countably infinite) P-framework.
	Suppose that $G$ has a cartesian $\Ck$-symmetric NAC-coloring $\delta$.
	Let $G'$ be the underlying unbraced subgraph of~$G$.
	
	Let $\bar{u}$ be any vertex of $G$.
	Let $R_1, R_2, \dots$ (resp.\ $B_1, B_2, \dots$) be the vertex sets
	of the connected components of~$G_\red^\delta$ (resp.\ $G_\blue^\delta$)\footnote{Although we use ``\dots'', it might happen that there are only finitely many such components (for instance if $G$ is finite).},
	and let $\bar{R}$, resp.\ $\bar{B}$, be the \red, resp. \blue{} component containing $\bar{u}$.
	Note that the vertex sets of these components are identical in~$G$ and $G'$.
	We define a map $\rho_\red:\{R_1, R_2, \dots\}\rightarrow\RR^2$ as follows:
	for $R_\ell$, let $W$ be any walk in the unbraced graph $G'$ from $\bar{u}$ to a vertex in $R_\ell$ and
	\begin{equation*}
		\rho_\red (R_\ell) = \sum_{\substack{\oriented{w_1}{w_2}\in W \\ \delta(w_1w_2)=\blue }} (\rho(w_2)-\rho(w_1))\,.
	\end{equation*}
	Notice that \cite[Lemma~4.3]{Grasegger2020} (or its infinite analogue) guarantees that $\rho_\red$ is well-defined,
	namely, the sum is independent of the choice of $W$ and the vertex in $R_\ell$.
	We define $\rho_\blue:\{B_1, \dots, B_n\}\rightarrow\RR^2$ analogously by swapping \red{} and \blue.
	
	For $t\in[0,2\pi]$ and $v\in V_G=V_{G'}$, where $v\in R_\ell\cap B_j$, let
	\begin{equation*}
		\rho_t(v) = \Rot(t)\cdot \rho_\red(R_\ell) + \rho_\blue(B_j)\,,
	\end{equation*}
	where $\Rot(t)$ is the clockwise rotation matrix by $t$ radians.
	The proof that $\rho_0(v)=\rho(v) - \rho(\bar{u})$ is the same as in \cite[Lemma~7.5]{Grasegger2020}
	as well as the fact that the lengths of edges in $E_{G}$
	are non-zero and constant along the flex.
	If $\bar{u}$ is an invariant vertex, then one can show that $\rho_t$ is $\Ck$-symmetric,
	but this is not the case in general.
	Hence, define the flex~$\widetilde{\rho}_t$ so that for each $v \in V_G$ we have
	\begin{align*}
		\widetilde{\rho}_t(v) &= \rho_t(v) - \frac{1}{k} \sum_{i=0}^{k-1} \rho_t(\omega^i\bar{u}) \\
			&= \Rot(t)\cdot \rho_\red(R_\ell) + \rho_\blue(B_j)
			- \frac{1}{k} \sum_{i=0}^{k-1}\left( \Rot(t)\cdot \rho_\red(\omega^i \bar{R}) + \rho_\blue(\omega^i \bar{B})\right)\,.
	\end{align*}
	Since the sum is independent of $v$, we have that the edge lengths are constant along the flex $\widetilde{\rho}_t$ as it is so for $\rho_t$.
	Recall that $\Rot_\omega$ is the clockwise rotation matrix $\Rot(2\pi/k)$.
	As
	\begin{align*}
		\frac{1}{k} \sum_{i=0}^{k-1} \rho_0(\omega^i\bar{u}) = \frac{1}{k} \sum_{i=0}^{k-1} \left(\rho(\omega^i\bar{u}) - \rho(\bar{u}) \right) = - \rho(\bar{u}) + \left(\frac{1}{k} \sum_{i=0}^{k-1} \Rot_\omega^i \right) \rho(\bar{u})  = - \rho(\bar{u}),
	\end{align*}
	and $\rho_t$ is a flex of the translation of $(G,\rho)$ by $\rho(\bar{u})$,
	we have that $\widetilde{\rho}_t$ is a flex of $\rho$.
	
	In order to prove that the flex $\widetilde{\rho}_t$ is $\Ck$-symmetric,
	we need the following identity (and analogous identity for $\rho_\blue$):
	\begin{align}
		\label{eq:ubarIdentity}
		k \cdot\rho_\red(\omega\bar{R}) = \sum_{i=0}^{k-1} \left(\rho_\red(\omega^i\bar{R}) - \Rot_\omega\cdot  \rho_\red(\omega^i\bar{R}) \right)\,.
	\end{align}
	The $0$-th summand is 0 since $\bar{u}\in\bar{R}$.
	If $W$ is a walk from $\bar{u}$ to $\omega \bar{u}$ in $G'$,
	then $\widetilde{W}_i = (W,\omega W, \dots, \omega^{i-1}W)$ is a walk from $\bar{u}$ to $\omega^i \bar{u}$ in $G'$.
	Hence, the $i$-th summand on the right hand side for $i\geq 1$ equals
	\begin{align*}
		&\sum_{\substack{\oriented{w_1}{w_2}\in \widetilde{W}_i \\ \delta(w_1w_2)=\blue }} (\rho(w_2)-\rho(w_1))
						  - \sum_{\substack{\oriented{w_1}{w_2}\in \widetilde{W}_i \\ \delta(w_1w_2)=\blue }}(\rho(\omega w_2)-\rho(\omega w_1)) \\
		&=	\sum_{\substack{\oriented{w_1}{w_2}\in W \\ \delta(w_1w_2)=\blue }}\sum_{j=0}^{i-1} 
							  \left( (\rho(\omega^j w_2)-\rho(\omega^j w_1)) - (\rho(\omega^{j+1} w_2)-\rho(\omega^{j+1} w_1)) \right)\\
		&= \sum_{\substack{\oriented{w_1}{w_2}\in W \\ \delta(w_1w_2)=\blue }}
							  \left( (\rho(w_2)-\rho(w_1)) - (\rho(\omega^{i} w_2)-\rho(\omega^{i} w_1))\right)\,.
	\end{align*}
	Notice that the last expression is zero for $i=0$.
	Therefore, the right hand side of \eqref{eq:ubarIdentity} is
	\begin{align*}
		\sum_{i=0}^{k-1} &\sum_{\substack{\oriented{w_1}{w_2}\in W \\ \delta(w_1w_2)=\blue }}
							  \left( (\rho(w_2)-\rho(w_1)) - (\rho(\omega^{i} w_2)-\rho(\omega^{i} w_1))\right)\\
				&=\sum_{i=0}^{k-1} \rho_\red(\omega \bar{R}) 
					- \sum_{\substack{\oriented{w_1}{w_2}\in W \\ \delta(w_1w_2)=\blue }}
					\underbrace{\sum_{i=0}^{k-1} (\rho(\omega^{i} w_2)-\rho(\omega^{i} w_1))}_{=0 \text{ due to the symmetry}}\,,
	\end{align*}
	which concludes the proof of the identity.
	Finally, for $v\in R_\ell \cap B_j$ and walks $W_1$ from $\bar{u}$ to $\omega \bar{u}$ and $W$ from $\bar{u}$ to $v$ in $G'$
	(hence $\omega W$ is from $\omega \bar{u}$ to $\omega v$) we have that $\rho_\red(\omega R_\ell)$ is equal to
	\begin{align*}
		\sum_{\substack{\oriented{w_1}{w_2}\in W_1 \\ \delta(w_1w_2)=\blue }} (\rho(w_2)-\rho(w_1))
		+\sum_{\substack{\oriented{w_1}{w_2}\in \omega W \\ \delta(w_1w_2)=\blue }} (\rho(w_2)-\rho(w_1)) = \rho_\red (\omega \bar{R}) + \Rot_\omega \rho_\red(R_\ell)\,,
	\end{align*}
	and similarly for $\rho_\blue(\omega B_j)$. As
	\begin{align*}
		 \Rot(t) &\left(\rho_\red(\omega R_\ell)- \frac{1}{k} \sum_{i=0}^{k-1}\rho_\red(\omega^i \bar{R})\right)
		  + \rho_\blue(\omega B_j) - \frac{1}{k} \sum_{i=0}^{k-1}\rho_\blue(\omega^i \bar{B}) \\
		  &=  \Rot(t) \left(\Rot_\omega \rho_\red(R_\ell) + \rho_\red (\omega \bar{R}) - \frac{1}{k} \sum_{i=0}^{k-1}\rho_\red(\omega^i \bar{R})\right) \\
		  &\qquad + \ \Rot_\omega \rho_\blue(B_j)+ \rho_\blue (\omega \bar{B}) - \frac{1}{k} \sum_{i=0}^{k-1}\rho_\blue(\omega^i \bar{B}) \\
		  &\stackrel{\text{\footnotesize{\eqref{eq:ubarIdentity}}}}{=}
		    \Rot(t) \left(\Rot_\omega \rho_\red(R_\ell)- \Rot_\omega\frac{1}{k} \sum_{i=0}^{k-1}\rho_\red(\omega^i \bar{R})\right)
		    + \ \Rot_\omega \rho_\blue(B_j)- \Rot_\omega\frac{1}{k} \sum_{i=0}^{k-1}\rho_\blue(\omega^i \bar{B}) \, ,
	\end{align*}
	we have  $\widetilde{\rho}_t(\omega v) = \Rot_\omega \widetilde{\rho}_t(v)$ as required.
\end{proof}

\begin{figure}[ht]
	\centering
	\begin{tikzpicture}
				\node[fvertex] (00) at (0,0) {};
				\node[fvertex] (10) at (1,0) {};
				\node[fvertex] (20) at (2,0) {};
				\node[fvertex] (01) at (0,1) {};
				\node[fvertex] (11) at (1,1) {};
				\node[fvertex] (21) at (2,1) {};
				\node[fvertex] (02) at (0,2) {};
				\node[fvertex] (12) at (1,2) {};
				\node[fvertex] (22) at (2,2) {};
				\draw[bedge] (00) to (10) (10) to (20) (20) to (21) (21) to (22)
							(22) to (12) (12) to (02) (02) to (01) (01) to (00);
				\draw[redge] (11) to (10) (11) to (21) (11) to (01) (11) to (12);
				\draw[edge,dashed] (-0.5,0.5) to (2.5,0.5) (-0.5,1.5) to (2.5,1.5) (0.5,-0.5) to (0.5,2.5) (1.5,-0.5) to (1.5,2.5);
	\end{tikzpicture}
	\caption{The $2 \times 2$ grid with its unique NAC-coloring (up to switching all red and blue edge colors). Its four ribbons are shown as gray dashed lines.}
	\label{fig:2by2symm}
\end{figure}
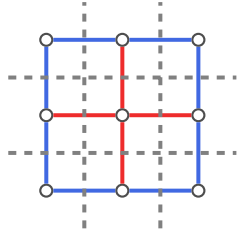

One aim of this section is to be able to determine whether a $\Ck$-symmetric braced  P-framework is $\Ck$-flexible solely by its bracing graph.
Our case example is the $2\times 2$ grid shown in \Cref{fig:2by2symm}.
While it does have a unique NAC-coloring (up to switching all red and blue edge colors) that is invariant under the group action of $\C_4$ as shown,
the NAC-coloring is neither cartesian nor $\C_4$-symmetric;
the two blue components, the central vertex and the outer cycle, are both invariant but are connected by four red edges.
It follows from \Cref{thm:bracedPframeworks,thm:finitesymm} that it cannot be $\C_4$-symmetrically flexible.
However the bracing graph is disconnected as it is just four isolated vertices.
As we shall soon see, this is actually not a problem,
for while the bracing graph of the P-framework is disconnected,
the quotient of the bracing graph by the symmetry of the ribbons is connected since it is a single vertex.
To capitalize on this observation, we define the following concept.

\begin{definition}
	Let $(G,\rho)$ be a $\Ck$-symmetric (countably infinite) braced P-framework with bracing graph~$\Gamma$.
	We define the \emph{quotient bracing graph} $\hat{\Gamma}$ to be the quotient graph of $\Gamma$ under the equivalence relation $r \sim r'$ if and only if $r' = \gamma r$ for some $\gamma \in \Ck$.
\end{definition}

With this definition we are finally ready for our main result of the section.
It is worth noting that this result is new for both finite and countably infinite graphs.

\begin{theorem}
	\label{thm:bracedCnPframeworks}
	Let $(G,\rho)$ be a $\Ck$-symmetric (countably infinite) braced P-framework.
	Then the following are equivalent:
	\begin{enumerate}[(i)]
		\item $(G,\rho)$ is $\Ck$-symmetrically rigid. \label{it:Cnrigid}
		\item $G$ has no cartesian $\Ck$-symmetric NAC-coloring. \label{it:CnnoNAC}
		\item The quotient bracing graph of $G$ is connected. \label{it:Cnconnected}
	\end{enumerate}
\end{theorem}

\begin{proof}
	$\neg \ref{it:CnnoNAC} \Rightarrow \neg \ref{it:Cnrigid}$:
	The statement holds by \Cref{lem:CnCartNACconstruction}.
	
	$\neg \ref{it:Cnrigid} \Rightarrow \neg \ref{it:CnnoNAC}$:
	First suppose that $(G,\rho)$ is finite.
	By \Cref{lem:active},
	there exists an algebraic motion $\C$ w.r.t.~the restricting set given by \cref{eq:parallelograms,eq:symmetric}, and an active NAC-coloring $\delta$ w.r.t.~$\C$.
	It was shown in \cite[Lemma 1]{RotSymmetry} that $\delta$ is $\Ck$-symmetric,
	and it follows immediately from \cref{eq:parallelograms} that $\delta$ must also be cartesian.
	
	Now suppose $(G,\rho)$ is a flexible infinite $\Ck$-symmetric P-framework.
	Let $\alpha \colon {[0,1)} \rightarrow (\RR^2)^{V_G}$ be the corresponding non-trivial $\Ck$-symmetric flex.
	Since the flex is non-trivial, there exist edges $e_1$ and $e_2$
	such that the angle between these two edges changes along the flex.
	Let $(G_n)_{n \in \NN}$ be a $\Ck$-symmetric subgraph tower of $G$ such that
	$e_1,e_2 \in E_{G_n}$.
	For any $n \in \NN$, projecting the flex~$\alpha$ to $G_n$ 
	gives a $\Ck$-symmetric flex $\alpha_n$ of~$G_{n}$.
	By a similar method to that used in \Cref{thm:bracedPframeworks},
	we may choose $\varepsilon_n \in (0,1)$ for every $n \in \NN$ so that $\alpha_n(t)$ is a $\Ck$-symmetric P-framework for every $t \in [0,\varepsilon_n)$.
	For each $n \in \NN$,
	apply \Cref{lem:active} to the restriction of $\alpha_n$ to $[0,\varepsilon_n)$ and the restricting set given by \cref{eq:parallelograms,eq:symmetric} to obtain a NAC-coloring $\delta_n$ with $\delta(e_1) \neq \delta(e_2)$.
	Each $\delta_n$ is $\Ck$-symmetric by \cite[Lemma 1]{RotSymmetry},
	and every ribbon of $G_n$ is monochromatic w.r.t.~$\delta_n$ because of \cref{eq:parallelograms}.
	We can now implement the methods showcased in \Cref{thm:bracedPframeworks} and \Cref{lem:CnNACfromSequence} with the sets
	\begin{align*}
		S_n = \{\delta\in\nac{G_n} : \delta \text{ is $\Ck$-symmetric} &\text{ with all ribbons being monochromatic} \\ 
					&\text{and } \delta(e_1)=\blue \neq \delta(e_2) \}
	\end{align*}
	to obtain the required cartesian $\Ck$-symmetric NAC-coloring of $G$.

	$\ref{it:Cnconnected} \Rightarrow \ref{it:CnnoNAC}$:	
	As the quotient bracing graph is connected,
	$\Gamma$ has at most $m$ connected components $\Gamma_1,\ldots,\Gamma_m$ for some $m|k$
	so that $\Gamma_i = \gamma_i \Gamma_1$ for some $\gamma_i \in \Ck$ for each $i \in \{1,\ldots,m\}$.
	By \Cref{lem:cartesianIffRibbonsMonochromatic} all edges in the same ribbon have the same color,
	hence so do all the edges of ribbons in the same connected component of the bracing graph.
	Because of $\Ck$ symmetry, the $\Ck$-symmetric NAC-coloring must be invariant under the group symmetry.
	This implies every edge has the same color, contradicting that $G$ has a cartesian $\Ck$-symmetric NAC-coloring.
	
	$\neg \ref{it:Cnconnected} \Rightarrow \neg \ref{it:CnnoNAC}$:
	Since the quotient NAC-coloring of $G$ is disconnected,
	there exists a red-blue vertex coloring $\delta_\Gamma$ of $\Gamma$ so that red (resp.~blue) vertices are only connected to red (resp.~blue) vertices,
	and $\delta_\Gamma(\gamma r) = \delta_\Gamma(r)$ for every $\gamma \in \Ck$.
	Define $\delta$ to be the red-blue edge coloring where for every edge $e$ in a ribbon $r$ we have $\delta(e) = \delta_\Gamma(r)$.
	As noted in \cite[Theorem~4.5]{Grasegger2020},
	$\delta$ is a cartesian NAC-coloring of $G$,
	and it is immediate that $\delta(\gamma e) = \delta(e)$ for all $e \in E_G$ and $\gamma \in \Ck$.
	By \Cref{lem:PartInvCompsPframework},
	$\delta$ is also $\Ck$-symmetric as required.	
\end{proof}

\begin{figure}[ht]
\centering
\begin{tikzpicture}[scale=0.5,rotate=-18]
	\node[fvertex] (0) at (0.00000, 0.00000) {};
	\node[fvertex] (1) at (-1.6180, -4.7684e-7) {};
	\node[fvertex] (2) at (-0.80902, -0.58779) {};
	\node[fvertex] (3) at (-0.50000, -1.5388) {};
	\node[fvertex] (4) at (0.30902, -0.95106) {};
	\node[fvertex] (5) at (1.3090, -0.95106) {};
	\node[fvertex] (6) at (1.0000, 0.00000) {};
	\node[fvertex] (7) at (1.3090, 0.95106) {};
	\node[fvertex] (8) at (0.30902, 0.95106) {};
	\node[fvertex] (9) at (-0.50000, 1.5388) {};
	\node[fvertex] (10) at (-0.80902, 0.58779) {};
	\node[fvertex] (11) at (-1.3090, -0.95106) {};
	\node[fvertex] (12) at (0.50000, -1.5388) {};
	\node[fvertex] (13) at (1.6180, 1.9073e-6) {};
	\node[fvertex] (14) at (0.50000, 1.5388) {};
	\node[fvertex] (15) at (-1.3090, 0.95106) {};
	\node[fvertex] (16) at (-2.6180, -9.5367e-7) {};
	\node[fvertex] (17) at (-2.3090, 0.95105) {};
	\node[fvertex] (18) at (-0.80902, -2.4899) {};
	\node[fvertex] (19) at (-1.6180, -1.9021) {};
	\node[fvertex] (20) at (2.1180, -1.5388) {};
	\node[fvertex] (21) at (1.3090, -2.1266) {};
	\node[fvertex] (22) at (2.1180, 1.5388) {};
	\node[fvertex] (23) at (2.4270, 0.58779) {};
	\node[fvertex] (24) at (-0.80902, 2.4899) {};
	\node[fvertex] (25) at (0.19098, 2.4899) {};
	\node[fvertex] (26) at (-2.3090, -0.95106) {};
	\node[fvertex] (27) at (0.19098, -2.4899) {};
	\node[fvertex] (28) at (2.4271, -0.58778) {};
	\node[fvertex] (29) at (1.3090, 2.1266) {};
	\node[fvertex] (30) at (-1.6180, 1.9021) {};
	\node[fvertex] (31) at (-2.6180, -1.9021) {};
	\node[fvertex] (32) at (1.0000, -3.0777) {};
	\node[fvertex] (33) at (3.2361, 1.9073e-6) {};
	\node[fvertex] (34) at (1.0000, 3.0777) {};
	\node[fvertex] (35) at (-2.6180, 1.9021) {};
	\node[fvertex] (36) at (-1.8090, -2.4899) {};
	\node[fvertex] (37) at (1.9073e-6, -3.0777) {};
	\node[fvertex] (38) at (1.8090, -2.4899) {};
	\node[fvertex] (39) at (2.9270, -0.95105) {};
	\node[fvertex] (40) at (2.9270, 0.95106) {};
	\node[fvertex] (41) at (1.8090, 2.4899) {};
	\node[fvertex] (42) at (-1.9073e-6, 3.0777) {};
	\node[fvertex] (43) at (-1.8090, 2.4899) {};
	\node[fvertex] (44) at (-2.9271, 0.95105) {};
	\node[fvertex] (45) at (-2.9270, -0.95106) {};
	\node[fvertex] (46) at (-3.4271, -0.58779) {};
	\node[fvertex] (47) at (-3.7361, -1.5388) {};
	\node[fvertex] (48) at (-0.49999, -3.4410) {};
	\node[fvertex] (49) at (0.30902, -4.0287) {};
	\node[fvertex] (50) at (3.1180, -1.5388) {};
	\node[fvertex] (51) at (3.9271, -0.95105) {};
	\node[fvertex] (52) at (2.4270, 2.4899) {};
	\node[fvertex] (53) at (2.1180, 3.4410) {};
	\node[fvertex] (54) at (-1.6180, 3.0777) {};
	\node[fvertex] (55) at (-2.6180, 3.0777) {};
	\node[fvertex] (56) at (-3.4270, -2.4899) {};
	\node[fvertex] (57) at (1.3090, -4.0287) {};
	\node[fvertex] (58) at (4.2361, 0.00000) {};
	\node[fvertex] (59) at (1.3090, 4.0287) {};
	\node[fvertex] (60) at (-3.4271, 2.4899) {};
	\node[fvertex] (61) at (-2.6180, -3.0777) {};
	\node[fvertex] (62) at (2.1180, -3.4410) {};
	\node[fvertex] (63) at (3.9270, 0.95106) {};
	\node[fvertex] (64) at (0.30902, 4.0287) {};
	\node[fvertex] (65) at (-3.7361, 1.5388) {};
	\node[fvertex] (66) at (-1.6180, -3.0777) {};
	\node[fvertex] (67) at (2.4271, -2.4899) {};
	\node[fvertex] (68) at (3.1180, 1.5388) {};
	\node[fvertex] (69) at (-0.50000, 3.4410) {};
	\node[fvertex] (70) at (-3.4270, 0.58778) {};
	\node[fvertex] (71) at (-1.3090, -4.0287) {};
	\node[fvertex] (72) at (3.4270, -2.4899) {};
	\node[fvertex] (73) at (3.4271, 2.4899) {};
	\node[fvertex] (74) at (-1.3090, 4.0287) {};
	\node[fvertex] (75) at (-4.2361, -9.5367e-7) {};
	\node[fvertex] (76) at (-2.3090, -4.0287) {};
	\node[fvertex] (77) at (3.1180, -3.4410) {};
	\node[fvertex] (78) at (4.2361, 1.9021) {};
	\node[fvertex] (79) at (-0.50000, 4.6165) {};
	\node[fvertex] (80) at (-4.5451, 0.95105) {};
	\node[fvertex] (81) at (-0.50000, -4.6165) {};
	\node[fvertex] (82) at (4.2361, -1.9021) {};
	\node[fvertex] (83) at (3.1180, 3.4410) {};
	\node[fvertex] (84) at (-2.3090, 4.0287) {};
	\node[fvertex] (85) at (-4.5451, -0.95106) {};
	\node[fvertex] (86) at (-1.6180, -4.9798) {};
	\node[fvertex] (87) at (-0.80901, -5.5676) {};
	\node[fvertex] (88) at (4.2361, -3.0777) {};
	\node[fvertex] (89) at (5.0451, -2.4899) {};
	\node[fvertex] (90) at (4.2361, 3.0777) {};
	\node[fvertex] (91) at (3.9270, 4.0287) {};
	\node[fvertex] (92) at (-1.6180, 4.9798) {};
	\node[fvertex] (93) at (-2.6180, 4.9798) {};
	\node[fvertex] (94) at (-5.2361, -1.9073e-6) {};
	\node[fvertex] (95) at (-5.5451, -0.95106) {};
	\node[fvertex] (96) at (4.7684e-6, -4.9798) {};
	\node[fvertex] (97) at (4.7361, -1.5388) {};
	\node[fvertex] (98) at (2.9270, 4.0287) {};
	\node[fvertex] (99) at (-2.9271, 4.0287) {};
	\node[fvertex] (100) at (-4.7361, -1.5388) {};
	\node[fvertex] (101) at (1.0000, -4.9798) {};
	\node[fvertex] (102) at (5.0451, -0.58778) {};
	\node[fvertex] (103) at (2.1180, 4.6165) {};
	\node[fvertex] (104) at (-3.7361, 3.4409) {};
	\node[fvertex] (105) at (-4.4271, -2.4899) {};
	\node[fvertex] (106) at (-2.6180, -4.9798) {};
	\node[fvertex] (107) at (3.9271, -4.0287) {};
	\node[fvertex] (108) at (5.0451, 2.4899) {};
	\node[fvertex] (109) at (-0.80902, 5.5676) {};
	\node[fvertex] (110) at (-5.5451, 0.95105) {};
	\node[fvertex] (111) at (-2.9270, -4.0287) {};
	\node[fvertex] (112) at (2.9271, -4.0287) {};
	\node[fvertex] (113) at (4.7361, 1.5388) {};
	\node[fvertex] (114) at (-9.5367e-7, 4.9798) {};
	\node[fvertex] (115) at (-4.7361, 1.5388) {};
	\node[fvertex] (116) at (-3.7361, -3.4410) {};
	\node[fvertex] (117) at (2.1180, -4.6165) {};
	\node[fvertex] (118) at (5.0451, 0.58779) {};
	\node[fvertex] (119) at (1.0000, 4.9798) {};
	\node[fvertex] (120) at (-4.4271, 2.4899) {};
	\node[fvertex] (121) at (-1.8090, -5.5676) {};
	\node[fvertex] (122) at (4.7361, -3.4410) {};
	\node[fvertex] (123) at (4.7361, 3.4410) {};
	\node[fvertex] (124) at (-1.8090, 5.5676) {};
	\node[fvertex] (125) at (-5.8541, -9.5367e-7) {};
	\draw[redge](24)edge(54) (2)edge(3) (40)edge(63) (75)edge(46) (75)edge(70) (85)edge(47) (56)edge(47) (105)edge(100) (8)edge(7) 
	(14)edge(7) (52)edge(22) (64)edge(59) (74)edge(54) (13)edge(5) (24)edge(30) (60)edge(55) (58)edge(63) (89)edge(97) (1)edge(10) 
	(113)edge(108) (80)edge(65) (16)edge(26) (48)edge(18) (3)edge(4) (18)edge(27) (41)edge(53) (115)edge(110) (41)edge(34) (16)edge(17)
	 (34)edge(42) (81)edge(49) (65)edge(44) (0)edge(2) (22)edge(23) (65)edge(60) (116)edge(111) (18)edge(19) (16)edge(70) (48)edge(71)
	  (84)edge(55) (9)edge(14) (32)edge(37) (36)edge(61) (11)edge(3) (9)edge(15) (73)edge(68) (24)edge(69) (83)edge(53) (72)edge(67)
	   (32)edge(57) (0)edge(10) (6)edge(7) (12)edge(5) (66)edge(71) (40)edge(33) (50)edge(20) (96)edge(87) (20)edge(21) (68)edge(22) 
	   (120)edge(115) (56)edge(31) (9)edge(10) (1)edge(2) (1)edge(11) (33)edge(58) (8)edge(9) (114)edge(109) (49)edge(37) (64)edge(79)
	    (59)edge(53) (36)edge(31) (67)edge(20) (32)edge(38) (113)edge(118) (62)edge(38) (56)edge(61) (99)edge(93) (1)edge(15) (4)edge(5) 
	    (72)edge(50) (20)edge(28) (43)edge(55) (64)edge(42) (16)edge(46) (97)edge(102) (74)edge(69) (18)edge(66) (0)edge(6) (98)edge(91)
	     (34)edge(59) (51)edge(39) (57)edge(62) (0)edge(4) (49)edge(57) (73)edge(52) (3)edge(12) (24)edge(25) (5)edge(6) (114)edge(119)
	      (104)edge(99) (43)edge(35) (58)edge(51) (35)edge(60) (76)edge(61) (82)edge(51) (100)edge(95) (77)edge(62) (45)edge(31) 
	      (78)edge(63) (0)edge(8) (45)edge(47) (35)edge(44) (13)edge(7) (29)edge(22) (106)edge(111) (96)edge(101) (112)edge(107) 
	      (98)edge(103) (112)edge(117) (33)edge(39);
	\draw[bedge](29)edge(14) (66)edge(61) (30)edge(15) (12)edge(21) (68)edge(63) (19)edge(31) (81)edge(71) (109)edge(79) (86)edge(71)
	 (76)edge(71) (94)edge(95) (26)edge(11) (88)edge(107) (124)edge(93) (26)edge(31) (85)edge(95) (74)edge(92) (46)edge(47) (16)edge(44) 
	 (56)edge(116) (104)edge(60) (90)edge(108) (72)edge(82) (113)edge(63) (86)edge(87) (92)edge(93) (35)edge(30) (64)edge(114) (58)edge(118)
	  (124)edge(109) (122)edge(107) (18)edge(37) (25)edge(34) (73)edge(90) (24)edge(9) (106)edge(76) (125)edge(95) (57)edge(117) (107)edge(77)
	   (16)edge(1) (125)edge(110) (20)edge(38) (100)edge(47) (50)edge(51) (112)edge(62) (110)edge(94) (75)edge(94) (72)edge(88) (24)edge(43) 
	   (25)edge(14) (121)edge(106) (89)edge(82) (32)edge(21) (65)edge(70) (20)edge(5) (57)edge(101) (64)edge(69) (58)edge(102) (106)edge(86)
	    (99)edge(55) (19)edge(11) (52)edge(53) (40)edge(22) (54)edge(55) (120)edge(60) (75)edge(85) (98)edge(53) (67)edge(62) (32)edge(27)
	     (20)edge(39) (59)edge(119) (48)edge(49) (33)edge(28) (18)edge(36) (84)edge(93) (121)edge(87) (34)edge(29) (73)edge(78) (73)edge(83)
	      (41)edge(22) (88)edge(89) (33)edge(23) (17)edge(15) (123)edge(108) (22)edge(7) (96)edge(49) (18)edge(3) (80)edge(110) (80)edge(75)
	       (72)edge(77) (74)edge(79) (28)edge(13) (108)edge(78) (16)edge(45) (24)edge(42) (92)edge(109) (97)edge(51) (59)edge(103) (90)edge(91)
	        (123)edge(91) (13)edge(23) (17)edge(35) (27)edge(12) (89)edge(122) (81)edge(87) (74)edge(84) (91)edge(83) (61)edge(111) (56)edge(105)
	         (65)edge(115);

	\end{tikzpicture}
	\caption{A $\Ck$-symmetric cartesian NAC-coloring yielding the flex depicted in \Cref{fig:flexingPenrose}.}
	\label{fig:PenroseNAC}
\end{figure}
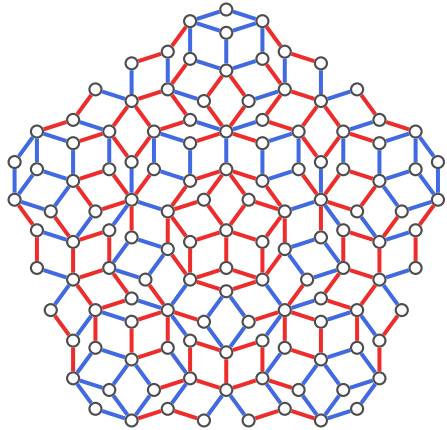

It is immediate that, similar to \Cref{thm:bracedPframeworks}, we can apply this result to either of the two rotationally symmetric Penrose frameworks.
In particular, a non-braced 5-fold symmetric Penrose framework will have a flex that preserves its symmetry; see \Cref{fig:PenroseNAC}.

\begin{proof}[Proof of \Cref{cor:SymmPenrose}]
	The quotient bracing graph of a Penrose framework with 5-fold symmetry is the countably infinite edgeless graph.
	The result now follows from \Cref{thm:bracedCnPframeworks}.
\end{proof}

\section{Flexibility of infinite complete bipartite frameworks}\label{sec:dixon}

\begin{figure}[ht]
	\centering
	\begin{tikzpicture}
		\draw[gridl, dashed] (-1.5,0)edge(2.3,0);
			\draw[gridl, dashed] (0,1.35)edge(0,-2.05);
			\node[fvertex] (2) at (1.8, 0) {};
			\node[fvertex] (5) at (-1., 0) {};
			\node[fvertex] (7) at (0.7, 0) {};
			\node[fvertex] (1) at (0, -1.75) {};
			\node[fvertex] (6) at (0,  1) {};
			\node[fvertex] (4) at (0, -0.8) {};
			\draw[edge]  (6)edge(5) (5)edge(4) (7)edge(4) (7)edge(6);
			\draw[edge] (2)edge(4) (2)edge(6);
			\draw[edge] (1)edge(5) (7)edge(1) ;
			\draw[edge] (2)edge(1);			
			\begin{scope}[xshift=3cm,yshift=-1.5cm]
				\draw[gridl, dashed] (1,0)edge(3,0) (1,3)edge(3,3) (0,1)edge(4,1) (0,2)edge(4,2);
				\draw[gridl, dashed] (0,1)edge(0,2) (4,1)edge(4,2) (1,0)edge(1,3) (3,0)edge(3,3); 
				\node[fvertex] (7) at (1,3) {};
				\node[fvertex] (5) at (4,1) {};
				\node[fvertex] (3) at (3,0) {};
				\node[fvertex] (1) at (3,3) {};
				\node[fvertex] (6) at (0,1) {};
				\node[fvertex] (8) at (0,2) {};
				\node[fvertex] (2) at (1,0) {};
				\node[fvertex] (4) at (4,2) {};
				\draw[edge] (5)edge(3) (2)edge(6) ;
				\draw[edge]  (2)edge(4)  (5)edge(7);
				\draw[edge] (2)edge(5) (7)edge(4) (3)edge(6);
				\draw[edge] (7)edge(6) (3)edge(4) ;
				\draw[edge]  (5)edge(1) (1)edge(6)  (1)edge(4) ;
				\draw[edge]  (1)edge(8) (2)edge(8) (3)edge(8) (7)edge(8);
			\end{scope}
	\end{tikzpicture}
	\caption{Dixon I (left) and II (right) constructions of flexible bipartite frameworks.}
	\label{fig:Dixon}
\end{figure}
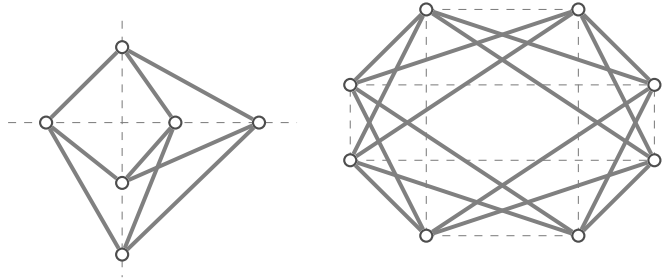

In \cite{Dixon}, Dixon showed two constructions of flexible frameworks for the complete bipartite graph $K_{4,4}$.
The first one (Dixon I) works for any finite bipartite graph and places the vertices of one part on $x$-axis
and the vertices of the other part on the $y$-axis; see \Cref{fig:Dixon}.
The second flexible realization possesses the dihedral symmetry of a rectangle, see also~\cite{Stachel}.
More than a hundred years later, it was proved by Maehara and Tokushige~\cite{Maehara2001}
that an injective realization of the complete bipartite graph $K_{m,n}$ vertices
with $m\geq 5$ and $n\geq 3$ is flexible if and only if the vertices are placed
as in the Dixon~I construction (modulo rigid motions).
Hence, Dixon~I construction is the only one that can apply in the case of flexible injective realizations of the countably infinite complete bipartite graph.
Contrary to the finite case, not all Dixon I realizations are flexible, as we will see in \Cref{thm:infiniteDixonI}.

Let us remark that for $K_{3,3}$, Walter and Husty \cite{WalterHusty} showed
that the realizations obtained by restricting the two Dixon's constructions
for $K_{4,4}$ are the only injective flexible ones,
see also \cite{GLSclassification} for an alternative proof.
However, there is also a third flexible injective realization on the sphere
besides the two analogues of Dixon's constructions~\cite{GGLSsphere}.

\begin{definition}
	Let $A = \{a_n : n \in \mathbb{N}\}$ and $B = \{b_n : n \in \mathbb{N}\}$ be vertex sets of the infinite complete bipartite graph $K_{A,B}$.
	Let $(K_{A,B},\rho)$ be the framework with $\rho(a_n)=(x_n,0)$ and $\rho(b_n)=(0,y_n)$, $\rho$ being an injective realization.
	We define $(K_{A,B},\rho)$ to be an \emph{infinite Dixon I linkage}.
\end{definition}

We now illustrate that, unlike with the finite case, there exist rigid infinite Dixon I linkages.

\begin{theorem}
	\label{thm:infiniteDixonI}
	Let $K_{A,B}$ be the complete bipartite graph with $A$ and $B$ countably infinite.
	Let $\rho$ be an injective realization of $K_{A,B}$.
	Then $(K_{A,B},\rho)$ is flexible if and only if it is congruent to an infinite Dixon I linkage
	with $\inf_{n\in \mathbb{N}} |x_n| \neq 0 \neq \inf_{n\in \mathbb{N}} |y_n|$.
\end{theorem}

\begin{proof}
	By the discussion above, a countably infinite complete bipartite framework can be flexible only
	if it is congruent to an infinite Dixon I linkage.
	We note that we may assume that $x_n ,y_n \geq 0$ for all $n \in \mathbb{N}$
	since this will equate to reflecting points through the $x$- and $y$-axis.
	We may also assume without loss of generality that $\inf_{n\in \mathbb{N}} x_n \geq \inf_{n\in \mathbb{N}} y_n$
	by reflecting the framework with respect to the line $y=x$.
	Hence our conditions to check are whether or not $\inf_{n\in \mathbb{N}} x_n$ is positive.
	
	First suppose that $c:=\inf_{n\in \mathbb{N}} x_n >0$.
	We show that Dixon~I construction works also for the infinite case under this assumption.
	For each $n \in \mathbb{N}$ and $t \in [0,1]$, the values
	$x_n(t):= \sqrt{x_n^2 - c^2(1-(1-t)^2)}$ and $y_n(t):= \sqrt{y_n^2 + c^2(1-(1-t)^2)}$
	are well-defined.
	With this, we define the continuous paths $\alpha_t(a_n) := (x_n(t),0)$ and $\alpha_t(b_n) := (0,y_n(t))$ for all $n \in \mathbb{N}$.
	For any $m,n \in \mathbb{N}$,
	the distance between $\alpha_t(a_m)$ and $\alpha_t(b_n)$ will remain constant as $t$ varies,
	hence $\alpha$ is a flex of $(K_{A,B},\rho)$.
	Since the vertices of $A$ and $B$ respectively are fixed to the axis,
	$\alpha$ is non-trivial as required.
	
	Now suppose that $\inf_{n\in \mathbb{N}} x_n = 0$.
	Let $(a_{n_j})_{j \in \mathbb{N}}$ and $(b_{n_k})_{k \in \mathbb{N}}$ be sequences
	where $x_{n_j} \rightarrow 0$ and  $y_{n_k} \rightarrow 0$ as $j,k \rightarrow \infty$.
	Suppose for contradiction that $\alpha$ is a non-trivial flex of $(K_{A,B},\rho)$.
	We first note that by the result of Maehara and Tokushige mentioned above, 
	we can assume that for all $n \in \mathbb{N}$ and $t \in [0,1]$,
	we have	$\alpha_t(a_n) = (x_n(t),0)$ and $\alpha_t(b_n) = (0,y_n(t))$
	for some continuous functions $x_n(t)$ and $y_n(t)$.
	For any $n \in \mathbb{N}$ we note that
	\begin{align*}
		\sup_{t \in [0,1]} x_n(t)^2 \leq \|\rho(a_n) - \rho(b_{n_k})\|^2 = x_n^2 + y_{n_k}^2 \rightarrow  x_n^2 \quad \text{ as } k \rightarrow \infty, \\
		\sup_{t \in [0,1]} y_n(t)^2 \leq \|\rho(a_{n_j}) - \rho(b_{n})\|^2 = x_{n_j}^2 + y_{n}^2 \rightarrow  y_n^2 \quad \text{ as } j \rightarrow \infty,
	\end{align*}
	thus $x_n(t) \leq x_n$ and $y_n(t) \leq y_n$ for all $t\in [0,1]$.
	If we choose any $m,n \in \mathbb{N}$ we see that for each $t \in [0,1]$ we have
	\begin{align*}
		x_m^2 +y_n^2=\|\rho(a_m) - \rho(b_n) \|^2 = \|\alpha_t(a_m) - \alpha_t(b_n) \|^2 = x_m(t)^2 + y_n(t)^2 \leq x_m^2 +y_n^2.
	\end{align*}
	It follows that $\alpha$ is constant, which contradicts non-triviality.
\end{proof}

We remark that all NAC-colorings of the complete countable infinite bipartite graph $K_{A,B}$
are such that any flexible realization constructed using \Cref{proposition:construction} will be non-injective.
This can be seen from the fact that this is the case already for all NAC-colorings
of the complete bipartite graph $K_{3,3}$, see \cite[Figure 10]{GLSclassification}.

\textsc{(SD) Johann Radon Institute for Computational and Applied Mathematics, Linz (RICAM)}\\

\textsc{(JL) Department of Applied Mathematics, Faculty of Information Technology, Czech Technical University in Prague}

\end{document}